\documentclass[12pt,titlepage]{article}
\usepackage{amssymb}
\input epsf
\newcommand{\proof}{{\noindent \bf Proof. }}
\newtheorem{thm}{Theorem}
\newtheorem{defi}{Definition}

\newtheorem{lem}{Lemma}[section]
\newtheorem{cor}[thm]{Corollary}
\newtheorem{prop}[thm]{Proposition}

\newcommand{\HH}{{\cal H}}

\def\2{\mathbb Z_2}
\def\ind{{\rm ind}}
\def\coind{{\rm coind}}
\def\susp{{\rm susp}}
\def\qed{$\Box$}

\newcommand\mbf[1]{\mbox{\boldmath$#1$}}
\newcommand\msbf[1]{\mbox{\boldmath\scriptsize$#1$}}

\title{Local chromatic number, Ky Fan's theorem, and circular colorings} 
\author{{\bf G\'abor Simonyi}\thanks{Research partially supported by the
Hungarian Foundation for Scientific Research Grant (OTKA) Nos.\ T037846 and
T046376.}   
$\qquad$\ \  {\bf G\'abor Tardos}\thanks{Research partially supported by the
Hungarian Foundation for Scientific Research Grant (OTKA) Nos.\ T037846 and
T046234.} \\ 
\medskip \\
Alfr\'ed R\'enyi Institute of Mathematics, \\
Hungarian Academy of Sciences, \\
1364 Budapest, POB 127, Hungary \\ 
\medskip \\
{\tt simonyi@renyi.hu} \ \ \  {\tt tardos@renyi.hu}}

\begin{document}
\maketitle

\begin{abstract}
The local chromatic number of a graph was introduced in \cite{EFHKRS}. It is
in between the chromatic and fractional chromatic numbers. This motivates the
study of the local chromatic number of graphs for which these quantities are
far apart. Such graphs include Kneser graphs, their vertex color-critical
subgraphs, the Schrijver (or stable Kneser) graphs; Mycielski graphs, and
their generalizations; and Borsuk graphs. We give more or less tight
bounds for the local chromatic number of many of these graphs.

We use an old topological result of Ky Fan \cite{kyfan} which generalizes the
Borsuk-Ulam theorem. It implies the
existence of a multicolored copy of the complete bipartite graph $K_{\lceil
t/2\rceil,\lfloor t/2\rfloor}$ in every proper coloring of many graphs whose
chromatic number $t$ is determined via a topological argument. 
(This was in particular noted for Kneser graphs by Ky Fan \cite{kyfan2}.) This
yields a lower bound of $\lceil t/2\rceil+1$ for the local chromatic number of
these graphs. We show this bound to be tight or almost tight in many cases.

As another consequence of the above we prove that the graphs considered here
have equal circular and ordinary chromatic numbers if the latter is even. This
partially proves a conjecture of Johnson, Holroyd, and Stahl and was
independently attained by F.~Meunier \cite{meunier}. We also show that odd
chromatic Schrijver graphs behave differently, their circular chromatic number
can be arbitrarily close to the other extreme. 
\end{abstract}

\section{Introduction}

The local chromatic number of a graph is defined in \cite{EFHKRS} as the
minimum number of colors that must appear within distance $1$ of a vertex. 
For the formal definition let $N(v)=N_G(v)$ denote the {\em neighborhood} of a
vertex $v$ in a graph $G$, that is, $N(v)$ is the set of vertices $v$ is
connected to.  

\begin{defi} \label{defi:lochr} {\rm(\cite{EFHKRS})}
The {\em local chromatic number} $\psi(G)$ of a graph $G$ is
$$\psi(G):=\min_c \max_{v\in V(G)} |\{c(u): u \in N(v)\}|+1,$$
where the minimum is taken over all proper colorings $c$ of $G$.
\end{defi}

The $+1$
term comes traditionally from considering ``closed neighborhoods''
$N(v)\cup\{v\}$ and results in a simpler form of the relations with other
coloring parameters. 

While the local chromatic number of a graph $G$ obviously cannot be more than
the chromatic number $\chi(G)$, somewhat surprisingly, it can be arbitrarily
less, cf.\ \cite{EFHKRS}, \cite{Fur}.
On the other hand, it was shown in
\cite{KPS} that
$$\psi(G)\ge \chi_f(G)$$
holds for any graph $G$, where
$\chi_f(G)$ denotes the fractional chromatic number of $G$. For the
definition and basic properties of the fractional chromatic number we refer to
the books \cite{SchU,GR}. 

This suggests to investigate the local chromatic
number of graphs for which the chromatic number and the fractional
chromatic number are far apart. This is our main goal in this paper. 

Prime examples of graphs with a large gap between the chromatic and the
fractional chromatic number are Kneser graphs and Mycielski graphs,
cf.\ \cite{SchU}. Other, closely related examples are provided by Schrijver
graphs, that are vertex color-critical induced subgraphs of Kneser graphs, and
many of the so-called generalized Mycielski graphs. In this introductory
section we focus on Kneser graphs and Schrijver graphs, Mycielski graphs and
generalized Mycielski graphs will be treated in detail in Subsection
\ref{subsect:gmyc}.

We recall that the Kneser graph $KG(n,k)$ is defined
for parameters $n\ge 2k$ as the graph with all $k$-subsets of an $n$-set as
vertices where two such vertices are connected if they represent disjoint
$k$-sets.  
It is a celebrated result of Lov\'asz \cite{LLKn} (see also
\cite{Bar,Gre}) proving the earlier conjecture of Kneser, that 
$\chi(KG(n,k))=n-2k+2$. For the fractional chromatic number one has
$\chi_f(KG(n,k))=n/k$ as easily follows from
the vertex-transitivity of $KG(n,k)$ via the Erd\H{o}s-Ko-Rado theorem,
see \cite{SchU,GR}.  

B\'ar\'any's proof \cite{Bar} of the Lov\'asz-Kneser theorem was generalized by
Schrijver \cite{Schr} who found a fascinating family of subgraphs of Kneser
graphs that are vertex-critical with
respect to the chromatic number. 

Let $[n]$ denote the set $\{1,2,\dots,n\}$.

\begin{defi} {\rm (\cite{Schr})}
The stable Kneser graph or {\em Schrijver graph} $SG(n,k)$ is defined as
follows.
\begin{eqnarray*}
V(SG(n,k))&=&\{A\subseteq [n]: |A|=k,\forall i:\ \{i,i+1\}\nsubseteq
A\ \ \hbox{\rm and}\ \ \{1,n\}\nsubseteq A\}\\
E(SG(n,k))&=&\{\{A,B\}: A\cap B=\emptyset\}
\end{eqnarray*}
\end{defi}

Thus $SG(n,k)$ is the subgraph induced by those vertices of $KG(n,k)$ that
contain no neighboring elements in the cyclically arranged basic set
$\{1,2,\dots,n\}$. These are sometimes called {\it stable $k$-subsets}.   
The result of Schrijver in \cite{Schr} is that
$\chi(SG(n,k))=n-2k+2(=\chi(KG(n,k))$, but
deleting any vertex of $SG(n,k)$ the chromatic number drops, i.e.,
$SG(n,k)$ is vertex-critical with respect to the chromatic number. 
Recently Talbot \cite{Tal} proved an Erd\H{o}s-Ko-Rado type result,  
conjectured by Holroyd and Johnson \cite{HJ}, 
which implies that the ratio of the number of vertices and
the independence number in $SG(n,k)$ is $n/k$. This gives $n/k\leq
\chi_f(SG(n,k))$ and equality follows by $\chi_f(SG(n,k))\leq
\chi_f(KG(n,k))=n/k$. Notice that $SG(n,k)$ is not vertex-transitive in
general. See more on Schrijver graphs in \cite{BjLo, LihLiu, Mat, Zie}. 

Concerning the local chromatic number it was 
observed by several people \cite{ZF,JK}, that $\psi(KG(n,k))\ge
n-3k+3$ holds, since the neighborhood of any vertex in $KG(n,k)$
induces a $KG(n-k,k)$ with chromatic number $n-3k+2$. Thus for $n/k$
fixed but larger than $3$, $\psi(G)$ goes to infinity with $n$ and $k$. In
fact, the results of \cite{EFHKRS} have a similar implication also for
$2<n/k\leq 3.$ Namely, it follows from those results, that if
a series of graphs $G_1,\dots, G_i,\dots$ is such that $\psi(G_i)$ is bounded, 
while $\chi(G_i)$ goes to infinity, then the number of colors to
be used in colorings attaining the local chromatic number grows at least 
doubly exponentially in the chromatic number. However, Kneser graphs
with $n/k$ fixed and $n$ (therefore also the chromatic number $n-2k+2$)
going to infinity cannot satisfy this, since the total number of vertices
grows simply exponentially in the chromatic number. 

The estimates mentioned in the previous paragraph are elementary. On the other
hand, all known proofs for $\chi(KG(n,k))\ge n-2k+2$ use topology or at least
have a topological flavor (see \cite{LLKn,Bar,Gre,MatCCA} to mention just a
few such proofs). They use (or at least, are inspired by) the Borsuk-Ulam
theorem. 

In this paper we use a stronger topological result due to Ky Fan \cite{kyfan}
to establish that all proper colorings of a $t$-chromatic Kneser, Schrijver
or generalized Mycielski graph contain a multicolored copy of a balanced
complete bipartite graph. This was noticed by Ky Fan for Kneser graphs
\cite{kyfan2}. We also show that the implied lower bound of $\lceil
t/2\rceil+1$ on the local chromatic number is tight or almost tight for many
Schrijver and generalized Mycielski graphs.

In the following section we summarize our main results in more detail.

\section{Results} \label{sect:results}

In this section we summarize our results
without introducing the topological notions needed to state the results in
their full generality. We will introduce the phrase that a graph $G$ is {\em
topologically $t$-chromatic} meaning that $\chi(G)\ge t$ and this fact can be
shown by a specific topological method, see Subsection \ref{subsect:bounds}.
Here we use this phrase only to emphasize the generality of the
corresponding statements, but the reader can
always substitute the phrase ``a topologically $t$-chromatic graph'' by ``a
$t$-chromatic Kneser graph'' or ``a $t$-chromatic Schrijver graph'' or by ``a
generalized Mycielski graph of chromatic number $t$''.
\medskip

Our general lower bound for the local chromatic number proven in
Section~\ref{sect:lowb} is the following.

\begin{thm} \label{thm:lowb}
If $G$ is topologically  $t$-chromatic for some $t\ge2$, then 
$$\psi(G)\ge\left\lceil t\over 2\right\rceil+1.$$
\end{thm}

This result on the local chromatic number is the immediate consequence of the 
Zig-zag theorem in Subsection \ref{subsect:kyfzag} that we state here in a
somewhat weaker form: 

\begin{thm} \label{thm:bip}
Let $G$ be a topologically $t$-chromatic graph and let $c$
be a proper coloring of $G$ with an arbitrary number of colors. Then there
exists a complete bipartite subgraph $K_{\lceil{t\over 2}\rceil,\lfloor{t\over
2}\rfloor}$ of $G$ all vertices of which receive a different color in $c$. 
\end{thm}

We use Ky Fan's generalization of the Borsuk-Ulam theorem \cite{kyfan} for the
proof.  The Zig-zag theorem was previously established for Kneser
graphs by Ky Fan \cite{kyfan2}.

We remark that J\'anos K\"orner \cite{JK} 
suggested to introduce a graph invariant $b(G)$ which is the size (number of
points) of the largest completely multicolored complete bipartite graph that
should appear in any proper coloring of graph $G$. It is obvious from the
definition that this parameter is bounded from above by $\chi(G)$ and bounded
from below by the local chromatic number $\psi(G)$. An obvious consequence of
Theorem~\ref{thm:bip} is that if $G$ is topologically $t$-chromatic, then
$b(G)\ge t$.

In Section~\ref{sect:upb} we show that Theorem~\ref{thm:lowb}
is essentially tight for several Schrijver and generalized Mycielski
graphs. In particular, this is always the case for a topologically
$t$-chromatic graph that has a {\em
wide} $t$-coloring as defined in Definition~\ref{defi:wide} in 
Subsection~\ref{ss:wide}.

As the first application of our result on wide colorings we show, that if the
chromatic number is fixed and odd, and the size of the Schrijver graph is
large enough, then Theorem~\ref{thm:lowb} is exactly tight:

\begin{thm} \label{thm:upb}
If $t=n-2k+2>2$ is odd and $n\ge4t^2-7t$ then
$$\psi(SG(n,k))=\left\lceil t\over 2\right\rceil+1.$$
\end{thm}

See Remark~4 in Subsection~\ref{subsect:schr} for a relaxed bound on $n$.  
The proof of Theorem~\ref{thm:upb} is combinatorial. 
It will also show that the claimed value of
$\psi(SG(n,k))$ can be attained with a coloring using $t+1$ colors and
avoiding the appearance of a totally multicolored $K_{\lceil{t\over
2}\rceil,\lceil{t\over 2}\rceil}.$ 
To appreciate the latter property, cf.\ Theorem~\ref{thm:bip}.  

Since $SG(n,k)$ is an induced subgraph of $SG(n+1,k)$
Theorem~\ref{thm:upb} immediately implies that 
for every fixed even $t=n-2k+2$ and $n, k$ large enough 
$$\psi(SG(n,k))\in\left\{{t\over 2}+1,{t\over 2}+2\right\}.$$ 

The lower bound for the local chromatic number in Theorem~\ref{thm:lowb} is
smaller than $t$ whenever $t\ge 4$ but Theorem~\ref{thm:upb} claims the
existence of Schrijver graphs with smaller local than ordinary
chromatic number only with chromatic number $5$ and up. In an upcoming paper
\cite{up} we prove that the local chromatic number of all $4$-chromatic
Kneser, Schrijver, or generalized Mycielski graphs is $4$. The reason is that
all these graphs satisfy a somewhat stronger property, they are {\em strongly}
topologically $4$-chromatic (see Definition \ref{defi:topres}). We will,
however, also show in \cite{up} that topologically $4$-chromatic graphs of
local chromatic number $3$ do exist. 

To demonstrate that requiring large $n$ and $k$ in Theorem~\ref{thm:upb} is
crucial we prove the following statement. 

\begin{prop} \label{prop:k2}
$\psi(SG(n,2))=n-2=\chi(SG(n,2))$ for every $n\ge 4$. 
\end{prop}

\medskip

As a second application of wide colorings we prove in
Subsection~\ref{subsect:gmyc} that Theorem~\ref{thm:lowb} 
is also tight for several generalized Mycielski graphs. These graphs will be
denoted by $M_{\msbf r}^{(d)}(K_2)$ where ${\mbf r}=(r_1,\dots,r_d)$ is a vector
of positive integers. See Subsection~\ref{subsect:gmyc} for the
definition. Informally, $d$ is the number of iterations and $r_i$ is the
number of ``levels'' in iteration $i$ of the generalized Mycielski
construction. $M_{\msbf r}^{(d)}(K_2)$ is proven to be $(d+2)$-chromatic
``because of a  topological reason'' by Stiebitz \cite{Stieb}. This
topological reason implies that these graphs are strongly topologically
$(d+2)$-chromatic. Thus Theorem~\ref{thm:lowb} applies and gives the lower
bound part of the following result.

\begin{thm} \label{thm:gmycspec7}
If ${\mbf r}=(r_1,\ldots,r_d)$, $d$ is odd, and $r_i\ge 7$ for all $i$, then 
$$\psi(M_{\msbf r}^{(d)}(K_2))=\left\lceil d\over2\right\rceil+2.$$   
\end{thm}

It will be shown in Theorem~\ref{thm:gmycspec4} that relaxing the $r_i\ge 7$
condition to $r_i\ge 4$ an only slightly weaker upper bound is still valid. 
As a counterpart we also show (see Proposition~\ref{prop:myc2} in
Subsection~\ref{subsect:gmyc}) that for the ordinary Mycielski construction,
which is the special case of ${\mbf r}=(2,\dots,2)$, the local
chromatic number behaves just like the chromatic number. 

The Borsuk-Ulam Theorem in topology is known to be equivalent (see Lov\'asz
\cite{LLgomb}) to the validity of a tight lower bound on the chromatic number
of graphs defined on the $n$-dimensional sphere, called Borsuk graphs. 
In Subsection~\ref{subsect:ctopI} we prove that the local chromatic number of
Borsuk graphs behaves similarly as that of the graphs already mentioned
above. In this subsection we also formulate a topological consequence of our
results on the tightness of the result of Ky Fan \cite{kyfan}. We also give a
direct proof for the same tightness result.

\medskip

The circular chromatic number $\chi_c(G)$ of a graph $G$ was introduced by
Vince \cite{Vin}, see Definition~\ref{defi:circ} in Section~\ref{sec:circ}. It
satisfies $\chi(G)-1<\chi_c(G)\le\chi(G)$. In
Section~\ref{sec:circ} we prove the following result using the Zig-zag
theorem.

\begin{thm}\label{thm:circ}
If $G$ is topologically $t$-chromatic and $t$ is even, then $\chi_c(G)\ge t$.
\end{thm}
\medskip

This theorem implies that $\chi_c(G)=\chi(G)$ if the chromatic number is even
for Kneser graphs, Schrijver graphs, generalized Mycielski graphs, and certain
Borsuk graphs. The result on Kneser and Schrijver graphs gives a partial
solution of a conjecture by Johnson, Holroyd, and Stahl \cite{JHS} and a
partial answer to a question of Hajiabolhassan and Zhu \cite{HZ}. These
results were independently obtained by Meunier \cite{meunier}. The result
on generalized Mycielski graphs answers a question of Chang, Huang, and Zhu
\cite{CHZ}. 

We will also discuss the circular chromatic number of odd chromatic Borsuk and 
Schrijver graphs showing that they can be close to one less than the chromatic
number. For generalized Mycielski graphs a similar result was proven by Lam,
Lin, Gu, and Song \cite{LLGS}, that we will also use.

\section{Lower bound} \label{sect:lowb}

\subsection{Topological preliminaries} \label{toppre}

The following is a brief overview of some of the topological concepts we
need. We refer to \cite{Bjhand,Hat} and \cite{Mat} for basic
concepts and also for a more detailed discussion of the notions and facts
given below.  
\smallskip
\par
\noindent
A {\em $\mathbb{Z}_2$-space} (or {\em involution space}) is a pair $(T,\nu)$
of a topological space $T$ and the involution $\nu:T\to T$, which is
continuous and satisfies that $\nu^2$ is the identity map. The points $x\in T$
and $\nu(x)$ are called {\em antipodal}. The involution
$\nu$ and the $\mathbb{Z}_2$-space $(T,\nu)$ are {\em free} if
$\nu(x)\ne x$ for all points $x$ of $T$. If the involution is understood
from the context we speak about $T$ rather than the pair $(T,\nu)$. This is
the case, in particular, for the unit sphere $S^d$ in ${\mathbb R}^{d+1}$
with the
involution given by the central reflection ${\mbf x}\mapsto-{\mbf x}$. 
A continuous map $f:S\to T$ between $\mathbb{Z}_2$-spaces $(S,\nu)$ and
$(T,\pi)$ is a {\em $\mathbb{Z}_2$-map} (or an {\em equivariant map}) if it
respects the respective
involutions, that is $f\circ\nu=\pi\circ f$. If such a map exists we write
$(S,\nu)\to(T,\pi)$. If
$(S,\nu)\to(T,\pi)$ does not hold we write $(S,\nu)\not\to(T,\pi)$. If both
$S\to T$ and $T\to S$ we call the $\mathbb Z_2$-spaces $S$ and $T$ \ $\mathbb
Z_2$-equivalent and write $S\leftrightarrow T$.

We try to avoid using homotopy equivalence and $\2$-homotopy equivalence 
(i.e., homotopy equivalence given by $\2$-maps), but we will have to use two
simple 
observations. First, if the $\2$-spaces $S$ and $T$ are $\2$-homotopy
equivalent, then $S\leftrightarrow T$. Second, if the space $S$ is
homotopy equivalent to a sphere $S^h$ (this relation is between topological
spaces, not $\2$-spaces), then for any involution $\nu$ we have
$S^h\to(S,\nu)$. 
\smallskip
\par
\noindent
The $\mathbb{Z}_2$-index of a $\mathbb{Z}_2$-space $(T,\nu)$ is defined 
(see e.g.\ \cite{MZ,Mat}) as 
$${\rm ind}(T,\nu):=\min\{d\ge 0:(T,\nu)\to S^d\},$$
where ${\rm ind}(T,\nu)$ is set to be $\infty$ if
$(T,\nu)\not\to S^d$ for all $d$. 
\smallskip
\par
\noindent
The $\mathbb{Z}_2$-coindex of a $\mathbb{Z}_2$-space $(T,\nu)$ is defined as
$${\rm coind}(T,\nu):=\max\{d\ge 0:S^d\to(T,\nu)\}.$$
If such a map exists for all $d$, then we set ${\rm coind}(T,\nu)=\infty$. 
Notice that if $(T,\nu)$ is not free, we have ${\rm
ind}(T,\nu)={\rm coind}(T,\nu)=\infty$. 
\medskip
\par
\noindent
Note that $S\to T$ implies $\ind(S)\le\ind(T)$ and $\coind(S)\le\coind(T)$. In
particular, $\2$-equivalent spaces have equal index and also equal coindex.

The celebrated Borsuk-Ulam Theorem can be stated in many equivalent
forms. Here we state three of them. For more equivalent versions and several
proofs we refer to \cite{Mat}. Here (i) and (ii) are standard forms of the
Borsuk-Ulam Theorem, while (iii) is clearly equivalent to (ii).
\medskip

\noindent
{\bf Borsuk-Ulam Theorem.}
{\em\begin{description}
\item[(i)] (Lyusternik-Schnirel'man version) Let $d\ge0$ and let $\HH$ be a
collection 
of open (or closed) sets covering $S^d$ with no $H\in\HH$ containing a pair of
antipodal points. Then $|\HH|\ge d+2$.
\item[(ii)] $S^{d+1}\not\to S^d$ for any $d\ge 0$.
\item[(iii)] For a $\mathbb{Z}_2$-space $T$ we have ${\rm ind}(T)\ge
{\rm coind}(T)$.  
\end{description}}
\medskip
\par
\noindent
The suspension $\susp(S)$ of a topological space $S$ is defined as the factor
of the space $S\times[-1,1]$ that identifies all the points in $S\times\{-1\}$
and identifies also the points in $S\times\{1\}$. If $S$ is a $\2$-space with
the involution $\nu$, then the suspension $\susp(S)$ is also a $\2$-space with
the involution $(x,t)\mapsto(\nu(x),-t)$. Any $\2$-map $f:S\to T$ naturally
extends to a $\2$-map $\susp(f):\susp(S)\to\susp(T)$ given by
$(x,t)\mapsto(f(x),t)$. We have $\susp(S^n)\cong S^{n+1}$ with a
$\2$-homeomorphism. These observations show the well known inequalities below.
\begin{lem}\label{suspension}
For any $\2$-space $S$ \ $\ind(\susp(S))\le\ind(S)+1$ and
$\coind(\susp(S))\ge\coind(S)+1$.
\end{lem}
\bigskip
\par
\noindent
A(n abstract) simplicial complex $K$ is a non-empty, hereditary set
system. That is, $F\in K$, $F'\subseteq F$ implies $F'\in K$ and we have
$\emptyset\in K$. In this paper we consider only finite simplicial complexes.
The non-empty sets in $K$ are called {\em simplices}.
We call the set
$V(K)=\{x:\{x\}\in K\}$
the set of {\em vertices} of $K$. In a {\em geometric realization} of $K$ a
vertex $x$ corresponds to a point $||x||$ in a Euclidean space, a simplex
$\sigma$ corresponds to its {\em body}, the convex hull of its vertices:
$||\sigma||={\rm conv}(\{||x||:x\in\sigma\})$. We assume that the points
$||x||$ for $x\in\sigma$ are affine independent, and so $||\sigma||$ is a
geometric simplex. We also assume that disjoint simplices have disjoint
bodies. The body of the complex $K$ is $||K||=\cup_{\sigma\in K}||\sigma||$, 
it is determined up to homeomorphism by $K$. Any point in $p\in||K||$ has
a unique representation as a convex combination $p=\sum_{x\in
V(K)}\alpha_x||x||$ such that $\{x:\alpha_x>0\}\in K$.

A map $f:V(K)\to V(L)$ is called simplicial if it maps simplices to
simplices, that is $\sigma\in K$ implies $f(\sigma)\in L$. In this case we
define $||f||:||K||\to||L||$ by setting $||f||(||x||)=||f(x)||$ for vertices
$x\in V(K)$ and taking an affine extension of this function to the
bodies of each of the simplices in $K$. If $||K||$ and $||L||$ are $\2$-spaces
(usually with an involution also given by simplicial maps), then we say that
$f$ is a {\em $\2$-map} if $||f||$ is a $\2$-map. If $||K||$ is a $\2$-space
we use $\ind(K)$ and $\coind(K)$ for $\ind(||K||)$ and $\coind(||K||)$,
respectively.

\smallskip
\par
\noindent
Following the papers \cite{AFL,Kriz,MZ} 
we introduce the {\it box complex} $B_0(G)$ for any finite graph $G$. See
\cite{MZ} for several similar complexes.
We define $B_0(G)$
to be a simplicial complex on the vertices $V(G)\times\{1,2\}$. For subsets
$S,T\subseteq V(G)$ we denote the set $S\times\{1\}\cup
T\times\{2\}$ by $S\uplus T$, following the convention of \cite{Mat,MZ}. 
For $v\in V(G)$ we denote by $+v$ the vertex
$(v,1)\in\{v\}\uplus\emptyset$ and $-v$ denotes the vertex
$(v,2)\in\emptyset\uplus\{v\}$. We set
$S\uplus T\in B_0(G)$ if $S\cap T=\emptyset$ and the complete bipartite graph
with sides $S$ and $T$ is a subgraph of $G$. Note that $V(G)\uplus\emptyset$
and $\emptyset\uplus V(G)$ are simplices of $B_0(G)$.

The $\mathbb{Z}_2$-map $S\uplus T\mapsto T\uplus S$ acts simplicially on
$B_0(G)$. It makes the body of the complex a free $\mathbb{Z}_2$-space.

We define the {\em hom space} $H(G)$ of $G$ to be the subspace consisting of
those points $p\in||B_0(G)||$ that, when written as a convex combination
$p=\sum_{x\in V(B_0(G))}\alpha_x||x||$ with $\{x:\alpha_x>0\}\in B_0(G)$ give
$\sum_{x\in V(G)\uplus\emptyset}\alpha_x=1/2$. 

Notice that $H(G)$ can also be obtained as the body of a {\em cell complex}
$Hom(K_2,G)$, see \cite{BK}, or a simplicial complex $B_{chain}(G)$, see
\cite{MZ}. 

A useful connection between $B_0(G)$ and $H(G)$ follows from a combination of
results of Csorba \cite{Cs} and Matou\v{s}ek and Ziegler \cite{MZ}.
\begin{prop}\label{csorba}
$||B_0(G)||\leftrightarrow\susp(H(G))$
\end{prop}

\proof Csorba \cite{Cs}
proves the $\2$-homotopy equivalence of $||B_0(G)||$ and the suspension of the
body of yet another box complex $B(G)$ of $G$. As we mentioned, $\2$-homotopy
equivalence implies $\2$-equivalence. Matou\v{s}ek and Ziegler
\cite{MZ} prove the $\2$-equivalence of $||B(G)||$ and
$H(G)$. Finally for $\2$-spaces $S$ and $T$ if $S\to T$, then
$\susp(S)\to\susp(T)$, therefore $||B(G)||\leftrightarrow H(G)$ implies
$\susp(||B(G)||)\leftrightarrow\susp(H(G))$.
\hfill\qed

Note that Csorba \cite{Cs} proves, cf.\ also \v{Z}ivaljevi\'c \cite{Ziv},  the
$\2$-homotopy equivalence of $||B(G)||$ and $H(G)$, 
and therefore we could also claim $\2$-homotopy equivalence in
Proposition~\ref{csorba}.

\subsection{Some earlier topological bounds} \label{subsect:bounds}

A graph homomorphism is an edge preserving map from the vertex set of a graph 
$F$ to the vertex set of another graph $G$.
If there is a homomorphism $f$ from $F$ to $G$, then it
generates a simplicial map from $B_0(F)$ to $B_0(G)$ in the natural
way. This map is a $\mathbb{Z}_2$-map and thus it shows
$||B_0(F)||\to||B_0(G)||$. Here $||B_0(F)||\not\to||B_0(G)||$ can often be
proved using the indexes or coindexes of these complexes and it implies the
non-existence of a homomorphism from $F$ to $G$. A similar argument applies
with the spaces $H(\cdot)$ in place of $||B_0(\cdot)||$.  

Coloring a graph $G$ with $m$ colors can be considered as a graph homomorphism
from $G$ to the complete graph $K_m$. The box complex $B_0(K_m)$ is the
boundary complex of the $m$-dimensional {\em cross-polytope} (i.e., the
convex hull of the basis vectors and their negatives in ${\mathbb R}^m$), thus
$||B_0(K_m)||\cong S^{m-1}$ with a $\2$-homeomorphism and
$\coind(B_0(G))\le\ind(B_0(G))\leq m-1$ is necessary for $G$ being
$m$-colorable. Similarly, $\coind(H(G))\le{\rm ind}(H(G))\leq m-2$ is also
necessary for $\chi(G)\leq m$ since $H(K_m)$ can be obtained
from intersecting the boundary of the $m$-dimensional cross-polytope with the
hyperplane $\sum x_i=0$, and therefore $H(K_m)\cong S^{m-2}$ with a
$\2$-homeomorphism.  These four lower bounds on $\chi(G)$ can be arranged in a
single line of inequalities using Lemma~\ref{suspension} and
Proposition~\ref{csorba}:
\begin{equation}\label{eq:chib1}
\chi(G)\ge\ind(H(G))+2\ge\ind(B_0(G))+1\ge\coind(B_0(G))+1\ge\coind(H(G))+2 
\end{equation}

In fact, many of the known proofs of Kneser's conjecture can be interpreted as
a proof of an appropriate lower bound on the (co)index of one of the above
complexes.  
In particular, B\'ar\'any's simple proof \cite{Bar} exhibits a
map showing $S^{n-2k}\to H(KG(n,k))$ to conclude that ${\rm
coind}(H(KG(n,k)))\ge n-2k$ and thus $\chi(KG(n,k))\ge n-2k+2$.
The even simpler proof of Greene \cite{Gre} exhibits a map showing
$S^{n-2k+1}\to B_0(KG(n,k))$ to conclude that ${\rm coind}(B_0(KG(n,k)))\ge
n-2k+1$ and thus $\chi(KG(n,k))\ge n-2k+2$. Schrijver's proof \cite{Schr} of
$\chi(SG(n,k))\ge n-2k+2$ is a generalization of B\'ar\'any's and it also can
be interpreted as a proof of $S^{n-2k}\to H(SG(n,k))$.
We remark that the same kind of technique is used with other complexes
related to graphs, too. In particular, Lov\'asz's original proof \cite{LLKn}
can also be considered as exhibiting a $\mathbb{Z}_2$-map from $S^{n-2k}$ to
such a complex, different from the ones we consider
here. For a detailed discussion of several such complexes and their usefulness
in bounding the chromatic number we refer the reader to \cite{MZ}.
\medskip
\par
\noindent
The above discussion gives several possible ``topological reasons'' that can
force a graph to be at least $t$-chromatic. Here we single out two such
reasons. 
The statement of our results in Section~\ref{sect:results} becomes precise by
applying the conventions given by the following definition. 

\begin{defi} \label{defi:topres}
We say that a graph $G$ is {\em topologically $t$-chromatic} if
$${\rm coind}(B_0(G))\ge t-1.$$
We say that a graph $G$ is {\em strongly topologically $t$-chromatic} if
$${\rm coind}(H(G))\ge t-2.$$
\end{defi}

By Equation~(\ref{eq:chib1}) if a graph is strongly topologically
$t$-chromatic, then it is topologically $t$-chromatic, and if $G$ is
topologically $t$-chromatic, then $\chi(G)\ge t$. In an upcoming paper
\cite{up} we will show the existence of a graph for any $t\ge 4$ that is
topologically $t$-chromatic but not strongly topologically $t$-chromatic. We
will also show that the two notions have different consequences in terms of
the local chromatic number for $t=4$.

The notion that a graph is (strongly) topologically $t$-chromatic is useful,
as it applies to many widely studied classes of graphs. As we mentioned above,
B\'ar\'any \cite{Bar} and Schrijver \cite{Schr} establish this for
$t$-chromatic Kneser and Schrijver graphs. For the reader's convenience we
recall the proof here. See the analogous statement for
generalized Mycielski graphs and (certain finite subgraphs of the) Borsuk
graphs after we introduce those graphs.

\begin{prop}\label{b-s}{\rm(B\'ar\'any; Schrijver)} The $t$-chromatic Kneser
and Schrijver graphs are strongly topologically $t$-chromatic.
\end{prop}

\proof 
We need to prove that $SG(n,k)$ is strongly topologically
$(n-2k+2)$-chromatic, i.e., that ${\rm coind}(H(SG(n,k)))\ge n-2k$. The
statement for Kneser graphs follows. For ${\mbf x}\in S^{n-2k}$ let $H_{\msbf x}$
denote the open hemisphere in $S^{n-2k}$ around $\mbf x$. Consider an
arrangement of the elements of $[n]$ on $S^{n-2k}$ so that each open
hemisphere contains a stable $k$-subset, i.e., a vertex of $SG(n,k)$. It is
not hard to check that identifying $i\in[n]$ with ${\mbf v}_i/|{\mbf v}_i|$ for
${\mbf v}_i=(-1)^i(1,i,i^2,\dots,i^{n-2k})\in{\mathbb R}^{n-2k+1}$ provides
such an arrangement. For each
vertex $v$ of $SG(n,k)$ and ${\mbf x}\in S^{n-2k}$ let
$D_v({\mbf x})$ denote the smallest distance of a point in $v$ from the set
$S^{n-2k}\setminus H_{\msbf x}$ and let $D({\mbf x})=
\sum_{v\in V(SG(n,k))}D_v({\mbf x})$. Note that $D_v({\mbf x})>0$ if $v$ is
contained in $H_{\msbf x}$ and therefore $D({\mbf x})>0$ for all $\mbf x$. Let
$f({\mbf x}):={1\over2D({\msbf x})}\sum_{v\in V(SG(n,k))}D_v({\mbf x})||{+}v||+
{1\over2D(-{\msbf x})}\sum_{v\in V(SG(n,k))}D_v(-{\mbf x})||{-}v||$. This $f$ is a
${\mathbb Z}_2$-map $S^{n-2k}\to H(SG(n,k))$ proving the proposition. 
\hfill\qed

\subsection{Ky Fan's result on covers of spheres and the Zig-Zag theorem}
\label{subsect:kyfzag}

The following result of Ky Fan \cite{kyfan} implies the
Lyusternik-Schnirel'man version of the Borsuk-Ulam theorem. Here we state two
equivalent versions of the result, all in terms of sets covering the
sphere. See the original paper for another version generalizing another
standard form of the Borsuk-Ulam theorem.
\bigskip
\par
\noindent
{\bf Ky Fan's Theorem.}
{\em\begin{description}
\item[(i)]Let $\cal A$ be a system of open (or a finite system of closed)
subsets of $S^k$ covering the entire sphere. Assume a linear order $<$ is
given on $\cal A$ and all sets $A\in\cal A$ satisfy $A\cap-A=\emptyset$. Then
there are sets $A_1<A_2<\dots<A_{k+2}$ of $\cal A$ and a point ${\mbf x}\in
S^k$ such 
that $(-1)^i {\mbf x}\in A_i$ for all $i=1,\dots,k+2$.
\item[(ii)]Let $\cal A$ be a system of open (or a finite system of closed)
subsets of $S^k$ such that $\cup_{A\in\cal A}(A\cup-A)=S^k$. Assume a linear
order $<$ is given on $\cal A$ and all sets $A\in\cal A$ satisfy
$A\cap-A=\emptyset$. Then there are sets $A_1<A_2<\dots<A_{k+1}$ of $\cal A$
and a point ${\mbf x}\in S^k$ such that $(-1)^i {\mbf x}\in A_i$ for all
$i=1,\dots,k+1$. 
\end{description}}

The Borsuk-Ulam theorem is easily seen to be implied by version (i), that
shows in particular, that $|{\cal A}|\ge k+2$. 
We remark that \cite{kyfan} contains the above statements only about closed
sets. The statements on open sets can be deduced by a standard argument using
the compactness of the sphere. We also remark that version (ii) is formulated
a little differently in \cite{kyfan}. A place where one finds exactly the above
formulation (for closed sets, but for any $\2$-space) is Bacon's paper
\cite{Bacon}.   

\bigskip
\noindent
{\bf Zig-zag Theorem}
{\em Let $G$ be a topologically $t$-chromatic finite graph and let $c$ be an
arbitrary proper coloring of $G$ by an arbitrary number of colors. We assume
the colors are linearly ordered. Then $G$ contains a complete bipartite
subgraph $K_{\lceil{t\over 2}\rceil,\lfloor{t\over 2}\rfloor}$ such that $c$
assigns distinct colors to all $t$ vertices of this
subgraph and these colors appear alternating on the two sides of the bipartite
subgraph with respect to their order.} 
\medskip
\par
\noindent
\proof
We have $\coind(B_0(G))\ge t-1$, so there exists a $\2$-map $f:S^{t-1}\to
B_0(G)$. For any color $i$ we define a set $A_i\subset S^{t-1}$ letting ${\mbf
x}\in A_i$ if and only if for the minimal simplex $U_{\msbf x}\uplus V_{\msbf x}$
containing $f({\mbf x})$ there exists a vertex $z\in U_{\msbf x}$ with
$c(z)=i$. These sets are open, but they do not necessarily cover the entire
sphere 
$S^{t-1}$. Notice that $-A_i$ consists of the points ${\mbf x}\in S^{t-1}$ with
$-{\mbf x}\in A_i$, which happens if and only if there exists a vertex $z\in
U_{-\msbf x}$ with $c(z)=i$. Here $U_{-\msbf x}=V_{\msbf x}$. For every ${\mbf x}\in
S^{t-1}$ either $U_{\msbf x}$ or $V_{\msbf x}$ is not empty, therefore we have
$\cup_i(A_i\cup-A_i)=S^{t-1}$. Assume for a
contradiction that for a color $i$ we have $A_i\cap-A_i\ne\emptyset$ and let
$\mbf x$ be a point in the intersection. We have a vertex $z\in U_{\msbf x}$ and
a vertex $z'\in V_{\msbf x}$ with $c(z)=c(z')=i$. By the definition of
$B_0(G)$ the vertices $z$ and $z'$ are connected in $G$. This contradicts 
the choice of $c$ as a proper coloring. The contradiction shows that $A_i\cap
-A_i=\emptyset$ for all colors $i$.

Applying version (ii) of Ky Fan's theorem we get that for some colors
$i_1<i_2<\dots<i_t$ and a point ${\mbf x}\in S^{t-1}$ we have $(-1)^j{\mbf x}\in
A_{i_j}$ for $j=1,2,\dots t$. This implies the existence of vertices $z_j\in
U_{(-1)^j\msbf x}$ with $c(z_j)=i_j$. Now $U_{(-1)^j\msbf x}=U_{\msbf x}$ for even
$j$ and $U_{(-1)^j\msbf x}=V_{\msbf x}$ for odd $j$. Therefore the complete
bipartite graph with sides $\{z_j|\hbox{$j$ is even}\}$ and $\{z_j|\hbox{$j$
is odd}\}$ is a subgraph of $G$ with the required properties.
\hfill\qed
\medskip

This result was previously established for Kneser graphs in \cite{kyfan2}.
\medskip
\par
\noindent
{\em Remark 1.}
Since for any fixed coloring we are allowed to order the colors in an
arbitrary manner, the Zig-zag Theorem implies the existence of several totally
multicolored copies of $K_{\lceil{t\over 2}\rceil,\lfloor{t\over 2}\rfloor}$. 
For a uniform random order any fixed totally multicolored $K_{\lceil{t\over
2}\rceil,\lfloor{t\over 2}\rfloor}$ satisfies the zig-zag rule with
probability $1/{t\choose\lfloor t/2\rfloor}$ if $t$ is odd and with
probability $2/{t\choose t/2}$ if $t$ is even. Thus the Zig-zag Theorem
implies the existence of many different totally multicolored subgraphs
$K_{\lceil{t\over 2}\rceil,\lfloor{t\over 2}\rfloor}$ in $G$:
${t\choose\lfloor t/2\rfloor}$ copies for odd $t$ and ${t\choose t/2}/2$ copies
for even $t$.

In the computation above we do not consider two subgraphs different if they
are isomorphic with an isomorphism preserving the color of the vertices. With
this convention, if the coloring uses only $t$ colors we get a totally
multicolored $K_{\lceil{t\over 2}\rceil,\lfloor{t\over 2}\rfloor}$ subgraph
with all possible colorings, and the number of these different subgraphs is
exactly the lower bound stated.
\hfill$\Diamond$

\medskip
\par
\noindent
{\bf Proof of Theorems~\ref{thm:lowb} and \ref{thm:bip}.}

Theorems~\ref{thm:lowb} and \ref{thm:bip} are direct consequences of the
Zig-zag theorem. Indeed, any vertex of the $\lfloor t/2\rfloor$ side of the
multicolored complete bipartite graphs has at least $\lceil t/2\rceil$
differently colored neighbors on the other side. \hfill$\Box$
\medskip
\par
\noindent
{\em Remark 2.}
Theorem~\ref{thm:lowb} gives tight lower bounds for the local chromatic number
of topologically $t$-chromatic graphs for odd $t$ as several examples of the
next section will show. In the upcoming paper \cite{up} we will present
examples that show that the situation is similar for even values of
$t$. However, the graphs establishing this fact are {\em not} strongly
topologically $t$-chromatic, whereas the graphs showing tightness of
Theorem~\ref{thm:lowb} for odd $t$ are. This leaves open the question whether
$\psi(G)\ge t/2+2$ holds for all strongly topologically $t$-chromatic graphs
$G$ and even $t\ge4$. While we will prove this statement in \cite{up} for
$t=4$ we do not know the answer for higher values of $t$.
\hfill$\Diamond$

\section{Upper bound} \label{sect:upb}

In this section we present the combinatorial constructions that prove
Theorems~\ref{thm:upb} and \ref{thm:gmycspec7}. 
In both cases general observations on
wide colorings (to be defined below) prove useful. The upper bound in either
of Theorems~\ref{thm:upb} or \ref{thm:gmycspec7}
implies the existence of certain open covers of
spheres. These topological consequences and the local chromatic
number of Borsuk graphs are discussed in the last subsection of this section.

\subsection{Wide colorings}\label{ss:wide}

We start here with a general method to alter a $t$-coloring and get a
$(t+1)$-coloring showing that $\psi\le t/2+2$. It works if the original
coloring was wide as defined below.

\begin{defi}\label{defi:wide}
A vertex coloring of a graph is called {\em wide} if the end vertices
of all walks of length $5$ receive different colors.
\end{defi}

Note that any wide coloring is proper, furthermore any pair of vertices of
distance $3$ or $5$ receive distinct colors. Moreover, if a graph has a wide
coloring it does not contain a cycle of length $3$ or $5$. For graphs that do
not have cycles of length $3$, $5$, $7$, or $9$ any coloring is wide that
assigns different colors to vertices of distance $1$, $3$ or $5$
apart. Another equivalent definition (considered in \cite{GyJS}) is that a
proper coloring is wide if the neighborhood of any color class is an
independent set and so is the second neighborhood.

\begin{lem}\label{lem:5ut}
If a graph $G$ has a wide coloring using $t$ colors, then $\psi(G)\le\lfloor
t/2\rfloor+2$.
\end{lem}

\proof Let $c_0$ be the wide $t$-coloring of $G$. We alter this coloring by
switching the color of the neighbors of the troublesome vertices to a new
color. We define a vertex $x$ to be {\em troublesome} if $|c_0(N(x))|>t/2$.
Assume the color $\beta$ is not used in the coloring $c_0$. For $x\in V(G)$ we
let  
$$c(x)=\left\{\begin{array}{lll}\beta&&\hbox{if $x$ has a troublesome
neighbor}\\c_0(x)&&\hbox{otherwise.}\end{array}\right.$$

The color class $\beta$ in $c$ is the union of the neighborhoods of troublesome
vertices. To see that this is an independent set consider any two vertices $z$
and $z'$ of color $\beta$. 
Let $y$ be a troublesome neighbor of $z$ and let $y'$
be a troublesome neighbor of $z'$. Both $c_0(N(y))$ and $c_0(N(y'))$ contain
more than half of the $t$ colors in $c_0$, therefore these sets are not
disjoint. We have a neighbor $x$ of $y$ and a neighbor $x'$ of $y'$ satisfying
$c_0(x)=c_0(x')$. This shows that $z$ and $z'$ are not connected, as otherwise
the walk $xyzz'y'x'$ of length $5$ would have two end vertices in the same
color class.

All other color classes of $c$ are subsets of the corresponding
color classes in $c_0$, and are therefore independent. Thus $c$ is a proper
coloring. 

Any troublesome vertex $x$ has now all its neighbors recolored,
therefore $c(N(x))=\{\beta\}$. For the vertices of $G$ that are not troublesome
one has $|c_0(N(x))|\le t/2$ and $c(N(x))\subseteq c_0(N(x))\cup\{\beta\}$,
therefore $|c(N(x))|\le t/2+1$. Thus the coloring $c$
shows $\psi(G)\le t/2+2$ as claimed. \hfill$\Box$
\medskip

We note that the coloring $c$ found in the
proof uses $t+1$ colors and any vertex that sees the
maximal number $\lfloor t/2\rfloor+1$ of the colors in its neighborhood must
have a neighbor of color $\beta$. In particular, for odd $t$ one will always
find two vertices of the same color in any $K_{(t+1)/2,(t+1)/2}$ subgraph.

\subsection{Schrijver graphs} \label{subsect:schr}

In this subsection we prove Theorem~\ref{thm:upb} which shows that the local
chromatic number of Schrijver graphs with certain parameters are as low as
allowed by Theorem~\ref{thm:lowb}. We also prove Proposition~\ref{prop:k2}
to show that for other Schrijver graphs the local chromatic number agrees with
the chromatic number.

For the proof of Theorem \ref{thm:upb} we will use the following simple
lemma. 

\begin{lem} \label{partav}
Let $u,v\subseteq [n]$ be two vertices of $SG(n,k)$. If there is a walk of
length $2s$ between $u$ and $v$ in $SG(n,k)$ then $|v\setminus u|\leq s(t-2)$,
where $t=n-2k+2=\chi(SG(n,k))$. 
\end{lem}

\proof
Let $xyz$ be a length two walk in $SG(n,k)$. Since $y$ is disjoint from $x$, 
it contains all but at most $n-2k=t-2$ elements of $[n]\setminus x$. As $z$ is
disjoint from $y$ it can contain at most $t-2$ elements not contained in
$x$. This proves the statement for $s=1$. 

Now let $x_0x_1\dots x_{2s}$ be a $2s$-length walk between $u=x_0$ and
$v=x_{2s}$ and
assume the statement is true for $s-1$. Since $|v\setminus u|\leq |v\setminus
x_{2s-2}|+|x_{2s-2}\setminus u|\leq (t-2)+(s-1)(t-2)$, the proof is completed
by induction. 
\hfill$\Box$

\medskip

We remark that Lemma \ref{partav} remains true for $KG(n,k)$ with literally
the same proof, but we will need it for $SG(n,k)$, this is why it is stated
that way.  

\medskip

\noindent{\bf Theorem~\ref{thm:upb}} (restated) {\em
If $t=n-2k+2>2$ is odd and $n\ge4t^2-7t$, then
$$\psi(SG(n,k))=\left\lceil t\over 2\right\rceil+1.$$
} 

\proof We need to show that
$\psi(SG(n,k))=(t+3)/2$. Note that the $t=3$ case is trivial as all
$3$-chromatic 
graphs have local chromatic number $3$. The lower bound for the local
chromatic number follows from Theorem~\ref{thm:lowb} and
Proposition~\ref{b-s}. 

We define a wide coloring $c_0$ of $SG(n,k)$ using $t$
colors. From this Lemma~\ref{lem:5ut} gives the upper bound on $\psi(SG(n,k))$.

Let $[n]=\{1,\dots,n\}$ be partitioned into $t$ sets, each containing
an odd number of consecutive elements of
$[n]$. More formally, $[n]$ is partitioned into disjoint sets $A_1,\dots,A_t$,
where each $A_i$ contains consecutive elements 
and $|A_i|=2p_i-1$. We need $p_i\ge2t-3$ for the proof,
this is possible as long as $n\ge t(4t-7)$ as assumed.

Notice, that $\sum_{i=1}^t(p_i-1)=k-1$, and therefore any $k$-element subset
$x$ of $[n]$ must contain more than half (i.e., at least $p_i$) of the
elements in some $A_i$. We define our coloring $c_0$ by arbitrarily
choosing such an index $i$ as the color $c_0(x)$. This is a proper coloring
even for the graph $KG(n,k)$ since if two sets $x$ and $y$ both contain more
than half of the elements of $A_i$, then they are not disjoint. 

As a coloring of $KG(n,k)$ the coloring $c_0$ is not wide. We need to show
that the coloring $c_0$ becomes wide if we restrict it to the subgraph
$SG(n,k)$.

The main observation is the following: $A_i$ contains a single subset of
cardinality $p_i$ that does not contain two consecutive elements. Let $C_i$ be
this set consisting of the first, third, etc.\ elements of $A_i$. A vertex of
$SG(n,k)$ has no two consecutive elements, thus a vertex $x$ of $SG(n,k)$ of
color $i$ must contain $C_i$.

Consider a walk $x_0x_1\dots  x_5$ of length $5$ in
$SG(n,k)$ and let $i=c_0(x_0)$. Thus the set $x_0$ contains
$C_i$. By Lemma \ref{partav} $|x_4\setminus x_0|\leq 2(t-2)$. In particular,
$x_4$ contains all but at most $2t-4$ elements  of $C_i$. As $p_i=|C_i|\ge
2t-3$, this means $x_4\cap C_i\neq\emptyset.$ Thus the set $x_5$, which is
disjoint from $x_4$, cannot contain all elements of $C_i$, showing
$c_0(x_5)\neq i$. This proves that the coloring $c_0$
is wide, thus Lemma~\ref{lem:5ut} completes the proof
of the theorem. \hfill$\Box$
\medskip

Note that the smallest Schrijver graph for which the
above proof gives $\psi(SG(n,k))<\chi(SG(n,k))$ is $G=SG(65,31)$ with
$\chi(G)=5$ and $\psi(G)=4$. In Remark~4 below we show how the lower bound on
$n$ can be lowered somewhat. After that we show that some lower bound is
needed as $\psi(SG(n,2))=\chi(SG(n,2))$ for every $n$.

\medskip
\par
\noindent
{\em Remark 3.}
In \cite{EFHKRS} universal graphs $U(m,r)$ are defined for which it is
shown that a graph $G$ can be colored with $m$ colors such that the
neighborhood of every vertex contains fewer than $r$ colors if and only if a
homomorphism from $G$ to $U(m,r)$ exists.
The proof of Theorem~\ref{thm:upb} gives, for odd $t$, a
$(t+1)$-coloring of $SG(n,k)$ (for appropriately large $n$ and $k$ that give 
chromatic number $t$) for which no neighborhood contains more than
$(t+1)/2$ colors, thus establishing the existence of a homomorphism from
$SG(n,k)$ to $U(t+1,(t+3)/2)$. This, in particular, proves that
$\chi(U(t+1,(t+3)/2))\ge t$, which is a special case of Theorem~2.6 in
\cite{EFHKRS}. It is not hard to see that this inequality is actually an
equality. Further, by the composition of the appropriate maps, the
existence of this homomorphism also proves that $U(t+1,(t+3)/2)$ is strongly
topologically $t$-chromatic.  
\hfill$\Diamond$

\medskip
\par
\noindent
{\em Remark 4.}
For the price of letting the proof be a bit more complicated one can
improve upon the bound given on $n$ in Theorem~\ref{thm:upb}. In
particular, one can show that the same conclusion holds for odd $t$ and $n\ge
2t^2-4t+3$. More generally, we can show
$\psi(SG(n,k))\le\chi(SG(n,k))-m=n-2k+2-m$ provided that
$\chi(SG(n,k))\ge2m+3$ and $n\ge8m^2+16m+9$ or $\chi(SG(n,k))\ge4m+3$ and
$n\ge20m+9$. The smallest Schrijver graph for which we can prove that the
local chromatic number is smaller than the ordinary chromatic number is
$SG(33,15)$ with $\chi=5$ but $\psi=4$. It has $1496$ vertices. (In general,
one has $|V(SG(n,k))|={n\over k}{{n-k-1}\choose {k-1}}$, cf.\ Lemma 1 in
\cite{Tal}.) The smallest $n$ and
$k$ for which we can prove $\psi(SG(n,k))<\chi(SG(n,k))$ is for the graph
$SG(29,12)$ for which $\chi=7$ but $\psi\le6$.

We only sketch the proof. For a similar and more detailed proof see
Theorem~\ref{thm:gmycspec4}. The idea is again to take a basic coloring $c_0$
of $SG(n,k)$ and obtain a new coloring $c$ by recoloring to a new color some 
neighbors of those vertices $v$ for which $|c_0(N(v))|$ is too large. 
The novelty is that now we do not
recolor all such neighbors, just enough of them, and also the definition of
the basic coloring $c_0$ is a bit different. Partition $[n]$ into $t=n-2k+2$
intervals $A_1,\dots,A_t$, each of odd length as in the proof of
Theorem~\ref{thm:upb} and also define $C_i$ similarly to be
the unique largest subset of $A_i$ not containing consecutive
elements. For a vertex $x$ we define $c_0(x)$ to be the {\em smallest} $i$ for
which $C_i\subseteq x$. Note that such an $i$ must exist. Now we define when to
recolor a vertex to the new color
$\beta$ if our goal is to prove $\psi(SG(n,k))\leq b:=t-m$, where $m>0$.  
We let $c(y)=\beta$ iff $y$ is the neighbor of a vertex $x$ having at least
$b-2$ different colors {\em smaller} than $c_0(y)$ in its
neighborhood. Otherwise, $c(y)=c_0(y)$. It is clear that $|c(N(x))|\leq b-1$
is satisfied, the only problem we face is that $c$ may not be a proper
coloring. To avoid this problem we only need that the recolored vertices form
an independent set. 
For each vertex $v$ define the index set $I(v):=\{j: v\cap C_j=\emptyset\}$. 
If $y$ and $y'$ are recolored vertices then they are neighbors of
some $x$ and $x'$, respectively, where $I(x)$ contains $c_0(y)$ and at least
$b-2$ indices smaller than $c_0(y)$ and $I(x')$ contains $c_0(y')$ and at
least $b-2$ indices smaller than $c_0(y')$. 
Since $[n]-(x\cup y)=t-2$, there are at most $t-2$ elements in
$\cup_{j\in I(x)}C_j$ not contained in $y$. 
The definition of $c_0$ also implies that at least one element of $C_j$ is
missing from $y$ for every $j<c_0(y)$. 
Similarly, there are at most $t-2$
elements in $\cup_{j\in I(x')}C_j$ not contained in $y'$ and at least one
element of $C_j$ is missing from $y'$ for every $j<c_0(y').$ 
These conditions lead to $y\cap y'\neq \emptyset$ if the sizes
$|A_i|=2|C_i|-1$ are appropriately chosen. In particular, if $t\ge 2m+3$ and
$|A_t|\ge 1,\, |A_{t-1}|\ge 2m+3,\, |A_{t-2}|\ge \dots\ge |A_{t-(2m+2)}|\ge
4m+5$, or $t\ge 4m+3$  
and $|A_t|\ge 1,\, |A_{t-1}|\ge 3,\, |A_{t-2}|\ge \dots\ge |A_{t-(4m+2)}|\ge
5$,  
then the above argument leads to a proof of $\psi(SG(n,k))\leq t-m$. 
(It takes some further but simple argument why the last two intervals $A_i$
can be chosen smaller than the previous ones.)
These two possible choices of the interval sizes give the two general bounds
on $n$ we claimed sufficient for attaining $\psi(SG(n,k))\leq t-m$. The
strengthening of Theorem~\ref{thm:upb} is obtained by the $m=(t-3)/2$ special
case of the first bound.  
\hfill$\Diamond$

\bigskip

\noindent{\bf Proposition~\ref{prop:k2}} (restated) {\em
$\psi(SG(n,2))=n-2=\chi(SG(n,2))$ for every $n\ge 4$.
}
\medskip

\proof In the $n=4$ case $SG(n,2)$ consists of a single edge and the statement
of the proposition is trivial.
Assume for a contradiction that $\psi(SG(n,2))\le n-3$ for some $n\ge5$ and let
$c$ be a proper coloring of $SG(n,2)$ showing this with the minimal number of
colors. As $\chi(SG(n,2))=n-2$ and any coloring of a graph $G$ with exactly 
$\chi(G)$ colors cannot show $\psi(G)<\chi(G)$ the coloring $c$ uses at least
$n-1$ colors.

It is worth visualizing the vertices of $SG(n,2)$ as diagonals of an $n$-gon
(see \cite{BjLo}). In other words, $SG(n,2)$ is the complement of the line
graph of $D$, where $D$ is the  complement of the cycle $C_n$. The color
classes are independent sets in $SG(n,2)$, so they are either stars or
triangles in $D$.

We say that a vertex $x$ {\em sees} the color classes of its neighbors. By our
assumption every vertex sees at most $n-4$ color classes. 

Assume a color class consists of a single vertex $x$. As $x$ sees at most
$n-4$ of the at least 
$n-1$ color classes we can choose a different color for $x$. The
resulting coloring attains the same local chromatic number with fewer
colors. This contradicts the choice of $c$ and shows that no color class is a
singleton.

A triangle color class is seen by all other edges of $D$. 
A star color class with
center $i$ and at least three elements is seen by all vertices that, as edges
of $D$, are not incident to $i$. For star color classes of two edges there can
be one additional vertex not seeing the class. So every color class is seen by
all but at most $n-2$ vertices. We double count the pairs of a vertex $x$ and a
color class $C$ seen by $x$. On one hand every vertex sees at most $n-4$
classes. On the other hand all the color classes are seen by at
least $\left({n\choose2}-n\right)-(n-2)$ vertices. We have
$$(n-1)\left({n\choose2}-2n+2\right)\le\left({n\choose2}-n\right)(n-4),$$
and this contradicts our $n\ge5$ assumption. The contradiction proves the
statement.
\rightline{$\Box$}

\subsection{Generalized Mycielski graphs} \label{subsect:gmyc}

Another class of graphs for which the
chromatic number is known only via the topological method is formed by 
generalized Mycielski graphs, see \cite{GyJS,Mat,Stieb}.
They are interesting for us
also for another reason: there is a big gap between their fractional and
ordinary chromatic numbers (see \cite{LPU,Tar}), therefore the local
chromatic number can take its value from a large interval. 

Recall that the Mycielskian
$M(G)$ of a graph $G$ is the graph defined on
$(\{0,1\}\times V(G))\cup \{z\}$ with edge set $E(M(G))=\{\{(0,v),(i,w)\}:
\{v,w\}\in E(G), i\in \{0,1\}\}\cup \{\{(1,v),z\}: v\in V(G)\}$. Mycielski
\cite{Myc} used this construction to increase the chromatic number of a graph
while keeping the clique number fixed: $\chi(M(G))=\chi(G)+1$ and
$\omega(M(G))=\omega(G)$.

Following Tardif \cite{Tar}, the same construction can also be described as the
direct (also called categorical) product of $G$ with a path on three vertices
having a loop at one end and then identifying all vertices that have the other
end of the path as their first coordinate. Recall that the direct product of
$F$ and $G$ is a graph on $V(F)\times V(G)$ with an edge between $(u,v)$ and
$(u',v')$ if and only if
$\{u,u'\}\in E(F)$ and $\{v,v'\}\in E(G)$. The generalized Mycielskian of $G$
(called a cone over $G$ by Tardif \cite{Tar}) $M_r(G)$ is then defined by
taking the direct product of $P$ and $G$, where $P$ is a path on $r+1$ 
vertices having a loop at one end, and then identifying all the vertices in the
product with the loopless end of the path as their first coordinate. With this
notation $M(G)=M_2(G)$. These graphs were considered
by Stiebitz \cite{Stieb}, who proved that if $G$ is $k$-chromatic ``for a
topological reason'' 
then $M_r(G)$ is $(k+1)$-chromatic for a similar
reason. (Gy\'arf\'as, Jensen, and Stiebitz \cite{GyJS} also consider these
graphs and quote Stiebitz's argument a special case of which is also presented
in \cite{Mat}.)
The topological reason of Stiebitz is in different terms than
those we use in this paper but using results of \cite{BK} they imply
strong topological $(t+d)$-chromaticity for graphs obtained by $d$
iterations of the generalized Mychielski construction starting, e.g, from
$K_t$ or from a $t$-chromatic Schrijver graph.
More precisely, Stiebitz proved that the body of the so-called neighborhood
complex ${\cal N}(M_r(G))$ of $M_r(G)$, introduced in \cite{LLKn} by Lov\'asz, 
is homotopy equivalent to the suspension of $||{\cal N}(G)||$.
Since ${\rm susp}(S^n)\cong S^{n+1}$ this implies
that whenever $||{\cal N}(G)||$ is homotopy equivalent to an $n$-dimensional
sphere, then $||{\cal N}(M_r(G))||$ is homotopy equivalent to the
$(n+1)$-dimensional sphere. This happens, for example, if $G$ is a complete
graph, or an odd cycle. By a recent result of Bj\"orner and de~Longueville
\cite{BjLo} we also have a similar situation if $G$ is isomorphic to any
Schrijver graph $SG(n,k)$. Notice that the latter include complete graphs and
odd cycles.

It is known, that $||{\cal N}(F)||$ is homotopy
equivalent to $H(F)$ for every graph $F$, see Proposition 4.2 in \cite{BK}. 
All this implies that ${\rm coind}(H(M_r(G)))={\rm coind}(H(G))+1$ whenever
$H(G)$ is homotopy equivalent to a sphere, in particular, whenever $G$ is a
complete graph or an odd cycle, or, more generally, a Schrijver graph. 
It is very likely that Stiebitz's proof can be generalized to show that
$H(M_r(G))\leftrightarrow\susp(H(G))$ and therefore ${\rm
coind}(H(M_r(G)))\ge{\rm coind}(H(G))+1$ holds always. Here we restrict
attention to graphs $G$ with $H(G)$ homotopy equivalent to a sphere. 

For an integer vector ${\mbf r}=(r_1,\dots,r_d)$ with $r_i\ge1$ for all $i$
we let $M_{\msbf r}^{(d)}(G)=M_{r_d}(M_{r_{d-1}}(\dots M_{r_1}(G)\ldots))$
denote the graph obtained by a $d$-fold application of the generalized
Mycielski construction with respective parameters $r_1,\dots,r_d$.

\begin{prop} \label{prop:Stieb} {\rm (Stiebitz)} 
If $G$ is a graph for which $H(G)$ is homotopy
equivalent to a sphere $S^h$ with $h=\chi(G)-2$ (in particular, $G$ is a
complete graph or an odd cycle, or, more generally, a Schrijver graph)
and ${\mbf r}=(r_1,\dots,r_d)$ is arbitrary, then $M_{\msbf r}^{(d)}(G)$ is
strongly topologically $t$-chromatic for $t=\chi(M_{\msbf
r}^{(d)}(G))=\chi(G)+d$.
\hfill\qed
\end{prop}

\medskip

It is interesting to remark that $\chi(M_r(G))>\chi(G)$ does
not hold in general if $r\ge 3$, e.g., for $\overline C_7$, the complement
of the $7$-cycle, one has $\chi(M_3(\overline C_7))=\chi(\overline
C_7)=4$. Still, the result of Stiebitz implies that the sequence $\{\chi(M_{\msbf
r}^{(d)}(G))\}_{d=1}^{\infty}$ may avoid to increase only a finite number of
times. 
\medskip

The fractional chromatic number of Mycielski graphs were determined by Larsen,
Propp, and Ullman \cite{LPU}, who proved that 
$\chi_f(M(G))=\chi_f(G)+{1\over {\chi_f(G)}}$ 
holds for every $G$. This already shows that there is a large
gap between the chromatic and the fractional chromatic number of
$M_{\msbf r}^{(d)}(G)$ if $d$ is large enough and $r_i\ge 2$ for all $i$, since
obviously, $\chi_f(M_r(F))\leq \chi_f(M(F))$ holds if $r\ge 2$. The 
previous result was generalized by Tardif \cite{Tar} who showed that
$\chi_f(M_r(G))$ can also be expressed by $\chi_f(G)$ as 
$\chi_f(G)+{1\over{\sum_{i=0}^{r-1}(\chi_f(G)-1)^i}}$ whenever $G$ has at
least one edge. 
\medskip

First we show that for the original Mycielski construction the local chromatic
number behaves similarly to the chromatic number.

\begin{prop} \label{prop:myc2}
For any graph $G$ we have
$$\psi(M(G))=\psi(G)+1.$$
\end{prop}

\proof
We proceed similarly as one does in the proof of $\chi(M(G))=\chi(G)+1$.
Recall that $V(M(G))=\{0,1\}\times V(G)\cup\{z\}$.

For the upper bound consider a coloring $c'$ of $G$ establishing its local
chromatic number and let $\alpha$ and $\beta$ be two colors not used by $c'$. 
We define
$c((0,x))=c'(x)$, $c((1,x))=\alpha$ and $c(z)=\beta$. This proper coloring
shows $\psi(M(G))\le\psi(G)+1$.

For the lower bound consider an arbitrary proper coloring $c$ of $M(G)$. 
We have to show that some vertex must see at least $\psi(G)$ different
colors in its neighborhood. 

We define the coloring $c'$ of $G$ as follows:
$$c'(x)=\left\{\begin{array}{lll}c((0,x))&&\hbox{if }
c((0,x))\ne c(z)\\
c((1,x))&&\hbox{otherwise.}\end{array}\right.$$
It follows from the construction that $c'$ is a proper coloring of $G$. Note
that $c'$ does not use the color $c(z)$. 

By the definition of $\psi(G)$, there is some vertex $x$ of $G$ that has at
least $\psi(G)-1$ different colors in its neighborhood $N_G(x)$. If
$c'(y)=c(0,y)$ for all vertices $y\in N_G(x)$, then the vertex $(1,x)$ has all
these colors in its neighborhood, and also the additional color $c(z)$. If
however $c'(y)\ne c(0,y)$ for a neighbor $y$ of $x$, then the vertex $(0,x)$
sees all the colors $c'(N_G(x))$ in its neighborhood $N_{M(G)}(0,x)$, 
and also the additional
color $c(0,y)=c(z)$. In both cases a vertex has $\psi(G)$ different colors in
its neighborhood as claimed.
\hfill$\Box$
\medskip

We remark that $M_1(G)$ is simply the graph $G$ with a new vertex connected to
every vertex of $G$, therefore the following trivially holds.

\begin{prop} \label{prop:myc1}
For any graph $G$ we have
$$\psi(M_1(G))=\chi(G)+1.$$
\end{prop}
\hfill$\Box$
\medskip

For our first upper bound we apply Lemma~\ref{lem:5ut}. We use the following
result of Gy\'arf\'as, Jensen, and Stiebitz \cite{GyJS}. The lemma below is an
immediate generalization of the $l=2$ special case of Theorem~4.1 in
\cite{GyJS}. 
We reproduce the simple proof from \cite{GyJS} for the sake of completeness.

\begin{lem}\label{gyarfasek} {\rm (\cite{GyJS})}
If $G$ has a wide coloring with $t$ colors and $r\ge7$, then $M_r(G)$ has a
wide coloring with $t+1$ colors.
\end{lem}

\proof As there is a homomorphism from $M_r(G)$ to $M_7(G)$ if $r>7$ it is
enough to give the coloring for
$r=7$. We fix a wide $t$-coloring $c_0$ of $G$ and use the additional
color $\gamma$. The coloring of $M_7(G)$ is given as
$$c((v,x))=\left\{\begin{array}{lll}\gamma&&\hbox{$v$ is the vertex at
distance $3$, $5$ or $7$ from the loop}\\
c_0(x)&&\hbox{otherwise.}\end{array}\right.$$
It is straightforward to check that $c$ is a wide coloring.
\hfill$\Box$
\medskip

We can apply the results of Stiebitz and Gy\'arf\'as et al.\ recursively to
give tight or almost tight bounds for the local chromatic number of the graphs
$M_{\msbf r}^{(d)}(G)$ in many cases:

\begin{cor}\label{nagymyc}
If $G$ has a wide $t$-coloring and ${\mbf r}=(r_1,\ldots,r_d)$ with
$r_i\ge7$ for all $i$, then $\psi(M_{\msbf r}^{(d)}(G))\le\frac{t+d}2+2$.

If $H(G)$ is homotopy equivalent to a sphere $S^h$, then $\psi(M_{\msbf
r}^{(d)}(G))\ge\frac{h+d}2+2$.
\end{cor}

\proof For the first statement we apply Lemma~\ref{gyarfasek} recursively to
show that $M_{\msbf r}^{(d)}(G)$ has a wide $(t+d)$-coloring and then apply
Lemma~\ref{lem:5ut}.

For the second statement we apply the result of Stiebitz recursively to show
that $H(M_{\msbf r}^{(d)}(G))$ is homotopy equivalent to $S^{h+d}$. As noted in
the preliminaries this implies $\coind(H(M_{\msbf r}^{(d)}(G)))\ge h+d$. By
Theorem~\ref{thm:lowb} the statement follows.
\hfill$\Box$
\medskip

\noindent{\bf Theorem~\ref{thm:gmycspec7}} (restated)
{\em If ${\mbf r}=(r_1,\ldots,r_d)$, $d$ is odd, and $r_i\ge 7$ for all $i$,
then $$\psi(M_{\msbf r}^{(d)}(K_2))=\left\lceil d\over2\right\rceil+2.$$}   

\medskip
\par
\noindent
{\bf Proof.} 
Notice that for ${\mbf r}=(r_1,\ldots,r_d)$ with $d$ odd and $r_i\ge7$ for all
$i$ the lower and upper bounds of Corollary~\ref{nagymyc} give the exact value
for the local chromatic number $\psi(M_{\msbf r}^{(d)}(K_2))=(d+5)/2$. This
proves the theorem.
\hfill$\Box$
\smallskip
\par
\noindent
Notice that a similar argument gives the exact value of $\psi(G)$ for the
more complicated graph 
$G=M_{\msbf r}^{(d)}(SG(n,k))$ whenever $n+d$ is odd, $r_i\ge
7$ for all $i$, and $n\ge 4t^2-7t$ for $t=n-2k+2$. This follows from
Corollary~\ref{nagymyc} via the wide colorability of $SG(n,k)$ for $n\ge
4t^2-7t$ shown in the proof of Theorem~\ref{thm:upb} and Bj\"orner and
de~Longueville's result \cite{BjLo} about the homotopy equivalence of
$H(SG(n,k))$ to $S^{n-2k}$.  

We summarize our knowledge on $\psi(M_{\msbf
r}^{(d)}(K_2))$ after proving the following theorem, which shows that
almost the
same upper bound as in Corollary~\ref{nagymyc} is implied from the relaxed
condition $r_i\ge4$.  

\bigskip
\begin{thm} \label{thm:gmycspec4}
For ${\mbf r}=(r_1,\ldots,r_d)$ with $r_i\ge4$ for all $i$ one has 
$$\psi(M_{\msbf r}^{(d)}(G))\leq \psi(G)+\left\lfloor d\over 2\right\rfloor+2.$$
Moreover, for $G\cong K_2$, the following slightly sharper bound holds: 
$$\psi(M_{\msbf r}^{(d)}(K_2))\le\left\lceil d\over2\right\rceil+3.$$
\end{thm}
\medskip

\proof
We denote the vertices of $Y:=M_{\msbf r}^{(d)}(G)$ in accordance to the
description of the generalized Mycielski construction via graph products. 
That is, a vertex of $Y$ is a sequence $a_1a_2\dots a_du$ of length $(d+1)$,
where $\forall i: \ a_i\in \{0,1,\dots,r_i\}\cup \{*\}$, 
$u\in V(G)\cup \{*\}$ and if $a_i=r_i$ for some $i$ then necessarily $u=*$ and
$a_j=*$ for every $j>i$, and this is the only way $*$ can appear in a 
sequence. To define adjacency we denote by $\hat P_{r_i+1}$ the path on
$\{0,1,\dots,r_i\}$ where the edges are of the form $\{i-1,i\}, i\in
\{1,\dots, r_i\}$ and there is a loop at vertex $0$.  
Two vertices $a_1a_2\dots a_du$ and $a_1'a_2'\dots
a_d'u'$ are adjacent in $Y$ if and only if 
$$u=*\hbox{ or }u'=*\hbox{ or }\{u,u'\}\in E(G)\hbox{ and}$$
$$\forall i:\ \ a_i=*\hbox{ or }a_i'=*\hbox{ or }\{a_i, a_i'\}\in E(\hat
P_{r_i+1}).$$
\smallskip

\noindent
Our strategy is similar to that used in Remark~4. Namely, we give an original
coloring $c_0$ and identify 
the set of ``troublesome'' vertices for this coloring and recolor most of the
neighbors of these vertices to a new color.

Let us fix a coloring $c_G$ of $G$ with at most $\psi(G)-1$ colors in the
neighborhood of a vertex. Let the colors we use
in this coloring be called $0,-1,-2$, etc. Now we define $c_0$ as
follows.
$$c_0(a_1\dots a_du)=\left\{\begin{array}{lll}c_G(u)&&
\hbox{if }\forall i:a_i\leq 2\\
i&&\hbox{if $a_i\ge 3$ is odd and $a_j\leq 2$ for all $j<i$}\\
0&&\hbox{if $\exists i:a_i\ge4$ is even and $a_j\leq 2$ for all $j<i$}
\end{array}\right.$$
It is clear that vertices having the same color
form independent sets, i.e., $c_0$ is a proper coloring. 
Notice that if a vertex has neighbors of many different
``positive'' colors, then it must have many coordinates that are equal to $2$. 
Now we recolor most of the neighbors of these vertices. 

Let $\beta$ be a color not used by $c_0$ and set $c(a_1\dots a_du)=\beta$ if
$|\{i: a_i\hbox{ is odd}\}|> d/2$.
(In fact, it would be enough to give color $\beta$ only to those of the above
vertices, for which the first $\lfloor{d\over 2}\rfloor$ odd coordinates are
equal to $1$. We recolor more vertices for the sake of simplicity.) 
Otherwise, let $c(a_1\dots a_du)=c_0(a_1\dots a_du)$. 

First, we have to show that $c$ is proper. To this end we only have to
show that no pair of vertices getting color $\beta$ can be adjacent. If two
vertices, ${\mbf x}=x_1\dots x_dv_x$ and ${\mbf y}=y_1\dots y_dv_y$ are colored
$\beta$ then both have more than $d/2$ odd coordinates 
(among their first $d$ coordinates). Thus there is some common
coordinate $i$ for which $x_i$ and $y_i$ are both odd. This implies that they
cannot be adjacent. 

Now we show that for any vertex $\mbf a$ we have $|c(N({\mbf
a}))\cap\{1,\ldots,d\}|\le d/2$. Indeed, if $|c_0(N({\mbf
a}))\cap\{1,\ldots,d\}|>d/2$ we have ${\mbf a}=a_1\dots a_du$ with more than
$d/2$ coordinates $a_i$ that are even and
positive. Furthermore, the first $\lfloor d/2\rfloor$ of these
coordinates should be $2$. Let $I$ be the set of indices of these first
$\lfloor d/2\rfloor$ even and positive coordinates. We claim that $c(N({\mbf
a}))\cap\{1,\ldots,d\}\subseteq I$. This is so, since if a neighbor has an odd
coordinate somewhere outside $I$, then it cannot have $*$ at the positions of
$I$, therefore it has more than $d/2$ odd coordinates and it is recolored by
$c$ to the color $\beta$.  

It is also clear that no vertex can see more than $\psi(G)-1$ ``negative'' 
colors in its neighborhood in either coloring $c_0$ or $c$. Thus the
neighborhood of any vertex can contain at most $\lfloor
d/2\rfloor+(\psi(G)-1)+2$ colors, where the last $2$ is added because of the
possible appearance of colors $\beta$ and $0$ in the neighborhood. This proves
$\psi(Y)\le d/2+\psi(G)+2$ proving the first statement in the
theorem. 

For $G\cong K_2$ the above gives $\psi(M_{\msbf r}^{(d)}(K_2))\leq
\lfloor d/2\rfloor +4$ 
which implies the second statement for odd $d$. For even $d$ the bound of the
second statement is $1$ less. We can gain $1$ as follows. When defining $c$
let us recolor to $\beta$ those vertices ${\mbf a}=a_1\dots a_du$, too, for
which the number of odd coordinates $a_i$ is exactly ${d\over 2}$ and
$c_G(u)=-1$. The proof proceeds similarly as before but we gain $1$ by
observing that those vertices who see $-1$ can see only ${d\over 2}-1$
``positive'' colors.
\hfill$\Box$
\smallskip

We collect the implications of
Theorems~\ref{thm:gmycspec7}, \ref{thm:gmycspec4} and
Propositions~\ref{prop:myc2} and \ref{prop:myc1}. It would be
interesting to estimate the value $\psi(M_{\msbf r}^{(d)}(K_2))$ for the missing
case ${\mbf r}=(3,\ldots,3)$. We have $\lceil d/2\rceil+2\le\psi\le d+2$ in
this case.

\begin{cor} \label{cor:mycgap}
For ${\mbf r}=(r_1,\ldots,r_d)$ we have
$$\psi(M_{\msbf r}^{(d)}(K_2))=\left\{\begin{array}{lll}
(d+5)/2&&\hbox{if $d$ is odd and }\forall i:r_i\ge7\\
\lceil d/2\rceil+2\hbox{ or }\lceil d/2\rceil+3&&
\hbox{if }\forall i:r_i\ge4\\
d+2&&\hbox{if }r_d=1\hbox{ or }\forall i:r_i=2.
\end{array}\right.$$
\end{cor}
\hfill$\Box$
\medskip

\noindent
{\em Remark 5.} The improvement for even $d$ given in the last paragraph of the
proof of Theorem~\ref{thm:gmycspec4} can also be obtained in a different
way we explain here. Instead of changing the rule for recoloring, we can
enforce that a vertex can see only $\psi(G)-2$ negative colors. This can be
achieved by setting the starting graph $G$ to be $M_4(K_2)\cong C_9$ instead of
$K_2$ itself and coloring this $C_9$ with the pattern
$-1,0,-1,-2,0,-2,-3,0,-3$ along the cycle. One can readily check that every
vertex can see only one non-$0$ color in its neighborhood. 

The same trick can be used also if the starting graph is not $K_2$ or $C_9$,
but some large enough Schrijver graph of odd chromatic number. Coloring it as
in the proof of Lemma~\ref{lem:5ut} (using the wide coloring as given in the
proof of Theorem~\ref{thm:upb}), we arrive to the same phenomenon if we
use the new color $\beta=0$.
\hfill$\Diamond$

\medskip
\par
\noindent
{\em Remark 6.}
Gy\'arf\'as, Jensen, and Stiebitz \cite{GyJS} use generalized Mycielski graphs
to show that another graph they denote by $G_k$ is $k$-chromatic. The way they
prove it is that they exhibit a homomorphism from $M_{\msbf r}^{(k-2)}(K_2)$ to
$G_k$ for ${\mbf r}=(4,\dots,4)$. The existence of this homomorphism implies
that $G_k$ is strongly topologically $k$-chromatic, thus its 
local chromatic number is at least $k/2+1$. We do not know any non-trivial
upper bound for $\psi(G_k)$. Also note that \cite{GyJS} gives universal graphs
for the property of having a wide $t$-coloring. By Lemma~\ref{lem:5ut} this
graph has $\psi\le t/2+2$. On the other hand, since any graph with a wide
$t$-coloring admits a homomorphism to this graph, and we have seen the wide
colorability of some strongly topologically $t$-chromatic graphs, it is
strongly topologically $t$-chromatic, as well. This gives $\psi\ge t/2+1$.
\hfill$\Diamond$

\subsection{Borsuk graphs and the tightness of Ky Fan's theorem}
\label{subsect:ctopI}

\begin{defi} \label{defi:Bogr}
The Borsuk graph $B(n,\alpha)$ of parameters $n$ and $0<\alpha<2$ is the
infinite graph whose vertices are the points of the unit sphere in ${\mathbb
R}^n$ (i.e., $S^{n-1}$) and its edges connect the pairs of points
with distance at least $\alpha$. 
\end{defi}

One easily sees that $\chi(B(n,\alpha))\ge n+1$, and, as Lov\'asz
\cite{LLgomb} remarks, this statement is equivalent to the Borsuk-Ulam
theorem. For $\alpha\ge\sqrt{2+2/n}$ this lower bound is sharp, see
\cite{LLgomb, Mat} (cf.\ also the proof of Corollary \ref{cor:Borpsi} below).  

The local chromatic number of Borsuk graphs for large enough $\alpha$ can also
be determined by our methods. First we want to argue that Theorem
\ref{thm:lowb} is applicable for this infinite graph. Lov\'asz gives in
\cite{LLgomb} a finite graph $G_P\subseteq B(n,\alpha)$ which has the property
that its neighborhood complex ${\cal N}(G_P)$ is homotopy equivalent to
$S^{n-1}$. Now we can continue the argument the same way as in the previous
subsection: Proposition 4.2 in \cite{BK} states that ${\cal N}(F)$ is homotopy
equivalent to $H(F)$ for every graph $F$, thus ${\rm
coind}(H(G_P))\ge n-1$, i.e., $G_P$ is (strongly) topologically
$(n+1)$-chromatic. As $G_P\subseteq B(n,\alpha)$ we
have $\lceil{{n+3}\over 2}\rceil\leq\psi(G_P)\le\psi(B(n,\alpha))$ by Theorem
\ref{thm:lowb}. 

The following lemma shows the special role of Borsuk graphs among
strongly topologically t-chromatic graphs. It will also show that our earlier
upper bounds on the local chromatic number have direct implications for Borsuk
graphs. 

\begin{lem}\label{gh}
A finite graph $G$ is strongly topologically $(n+1)$-chromatic if and only if
for some  $\alpha<2$ there is a graph homomorphism from $B(n,\alpha)$ to $G$. 
\end{lem}

\proof
For the if part consider the finite graph $G_P\subseteq B(n,\alpha)$ given by
Lov\'asz \cite{LLgomb} satisfying $\coind(H(G_P))\ge n-1$.
If there is a homomorphism from $B(n,\alpha)$
to $G$, it clearly gives a homomorphism also from $G_P$ to $G$ which further
generates a $\2$-map from $H(G_P)$ to $H(G)$. This proves
$\coind(H(G))\ge n-1$. 

For the only if part, let $f:S^{n-1}\to H(G)$ be a $\2$-map. For a point
${\mbf x}\in S^{n-1}$ write $f({\mbf x})\in H(G)$ as the convex combination
$f({\mbf x})=\sum\alpha_v({\mbf x})||{+}v||+\sum\beta_v({\mbf x})||{-}v||$ of
the vertices of $||B_0(G)||$. Here the
summations are for the vertices $v$ of $G$, $\sum\alpha_v({\mbf
x})=\sum\beta_v({\mbf x})=1/2$, and $\{v:\alpha_v({\mbf
x})>0\}\uplus\{v:\beta_v({\mbf x})>0\}\in B_0(G)$. Note that $\alpha_v$ and
$\beta_v$ are continuous as $f$ is continuous and $\beta_v({\mbf
x})=\alpha_v(-{\mbf x})$ by the equivariance of $f$. Set
$\varepsilon=1/(2|V(G)|)$. 
For ${\mbf x}\in S^{n-1}$
select an arbitrary vertex $v=g({\mbf x})$ of $G$ with
$\alpha_v\ge\varepsilon$. We claim that $g$ is a graph homomorphism from
$B(n,\alpha)$ to $G$ if $\alpha$ is close enough to $2$. By compactness it is
enough to prove that if we have vertices $v$ and $w$ of $G$ and sequences ${\mbf
x}_i\to{\mbf x}$ and ${\mbf y}_i\to-{\mbf x}$ of points in $S^{n-1}$ with $g({\mbf
x}_i)=v$ and $g({\mbf y}_i)=w$ for all $i$, then $v$ and $w$ are connected in
$G$. But since $\alpha_v$ is continuous we have $\alpha_v({\mbf
x})\ge\varepsilon$ 
and similarly $\beta_w({\mbf x})=\alpha_w(-{\mbf x})\ge\varepsilon$ and so $+v$
and 
$-w$ are contained in the smallest simplex of $B_0(G)$ containing $f({\mbf x})$
proving that $v$ and $w$ are connected.
\hfill$\Box$

By Lemma \ref{gh} either of Theorems \ref{thm:upb} or \ref{thm:gmycspec7}
implies that the above given lower bound on $\psi(B(n,\alpha))$ is tight
whenever $\chi(B(n,\alpha))$ is odd, that is, $n$ is even, and $\alpha<2$ is
close enough to $2$. The following corollary uses Lemma~\ref{lem:5ut} directly
to have an explicit bound on $\alpha$. 

\begin{cor} \label{cor:Borpsi}
If $n$ is even and $\alpha_n\le\alpha<2$, then
$$\psi(B(n,\alpha))={n\over 2}+2,$$
where $\alpha_n=2\cos{\arccos(1/n)\over10}$.
\end{cor}

Note that $\alpha_n\le\alpha_2=2\cos(\pi/30)<1.99$.
\medskip

\proof The lower bound on $\psi(B(n,\alpha))$ follows from
the discussion preceding Lemma~\ref{gh}.
The upper bound follows from Lemma~\ref{lem:5ut} as long as we can
give a wide $(n+1)$-coloring $c_0$ of the graph $B(n,\alpha)$.

To this end we use the standard $(n+1)$-coloring of $B(n,\alpha)$
(see, e.g., \cite{LLgomb,Mat}).
Consider a regular simplex $R$ inscribed into the unit sphere $S^{n-1}$ and
color a point ${\mbf x}\in S^{n-1}$ by the facet of $R$ intersected by the
segment from the origin to $\mbf x$. If this segment meets a lower
dimensional face then we arbitrarily choose a facet containing this face.

We let $\varphi=2\arccos(\alpha/2)$. Clearly, $\mbf x$ and $\mbf y$ is connected
if and only if the length of the shortest arc on $S^{n-1}$ connecting $-\mbf x$
and $\mbf y$ is at most $\varphi$. Therefore $\mbf x$ and $\mbf y$ are
connected by a walk of length $5$ if and only if the length of this same
minimal arc is at most $5\varphi$. For the coloring $c_0$ the length of the
shortest arc between $-\mbf x$ and $\mbf y$ for two vertices $\mbf x$ and $\mbf y$
colored with the same color is at 
least $\arccos(1/n)$. Hence $c_0$ is wide. To make sure this holds for
$\alpha=\alpha_n$ one has to color the $n+1$ vertices of the simplex $R$ with
different colors. \hfill$\Box$
\medskip
\par
\noindent

Our investigations of the local chromatic number led us to consider the 
following function $Q(h)$. The question of its values was independently asked
by Micha Perles motivated by a related question of Matatyahu Rubin\footnote{We
thank Imre B\'ar\'any \cite{b} and Gil Kalai \cite{GK} for this information.}. 

\begin{defi}\label q
For a nonnegative integer parameter $h$ let $Q(h)$ denote the minimum $l$ for
which $S^h$ can be covered by open sets in such a way that no point of the
sphere is contained in more than $l$ of these sets and none of the covering
sets contains an antipodal pair of points.
\end{defi}

Ky Fan's theorem implies $Q(h)\ge {h\over 2}+1$.
Either of Theorems \ref{thm:upb} or \ref{thm:gmycspec7} implies 
the upper bound $Q(h)\leq {h\over 2}+2$. 
Using the concepts of Corollary~\ref{cor:Borpsi} and Lemma~\ref{lem:5ut} one
can give an explicit covering of the sphere $S^{2l-3}$ by open subsets
where no point is contained in more than $l$ of the sets and no set contains
an antipodal pair of points. In fact, the covering we give satisfies a stronger
requirement and proves that version (ii) of Ky Fan's theorem is 
tight, while version (i) is almost tight.

\begin{cor} \label{kyftight}
There is a configuration ${\cal A}$ of $k+2$ open (closed) sets such that
$\cup_{A\in\cal A}(A\cup-A)=S^k$, all sets $A\in\cal A$ satisfy
$A\cap-A=\emptyset$, and 
no ${\mbf x}\in S^k$ is contained in more
than $\left\lceil k+1\over2\right\rceil$ of these sets. 

Furthermore, for every $\mbf x$ the number of sets in $\cal A$ containing
either $\mbf x$ or $-\mbf x$ is at most $k+1$.
\end{cor}

\proof 
First we construct closed sets. Consider the unit sphere
$S^k$ in ${\mathbb R}^{k+1}$. Let $R$ be a regular simplex inscribed in the
sphere. Let $B_1,\ldots,B_{k+2}$ be the subsets of the sphere obtained by the
central projection of the facets of $R$. These closed sets cover
$S^k$. Let $C_0$ be the set of points covered by at least
$\left\lceil k+3\over2\right\rceil$ of the sets
$B_i$. Notice that $C_0$ is the union of the central projections of the
$\lfloor{{k-1}\over 2}\rfloor$-dimensional faces of $R$.  
For odd $k$ let $C=C_0$, while for even $k$ let $C=C_0\cup C_1$, where $C_1$
is the set of points in $B_1$ covered by exactly $k/2+1$ of the sets
$B_i$. Thus $C_1$ is the central projection of the ${k\over 2}$-dimensional
faces of a facet of $R$. Observe that $C\cap-C=\emptyset$. Take
$0<\delta<\hbox{dist}(C,-C)/2$ and let $D$ be the open $\delta$-neighborhood of
$C$ in $S^k$. For $1\le i\le k+2$ let $A_i=B_i\setminus D$. These closed sets
cover $S^k\setminus D$ and none of them contain a pair of antipodal points. As
$D\cap-D=\emptyset$ we have $\cup_{i=1}^{k+2}(A_i\cup-A_i)=S^k$. It is clear
that every point of the sphere is covered by at most $\left\lceil
k+1\over2\right\rceil$ of the sets $A_i$ proving the first statement of the
corollary.

For the second statement note that if each set $B_i$ contains at least one of
a pair of antipodal points, then one of these points
belongs to $C$ and is therefore not covered by any of the sets $A_i$. Note
also, that for odd $k$ the second statement follows also from the first.

To construct open sets as required we can simply take the open
$\varepsilon$-neighborhoods of $A_i$. For small enough $\varepsilon>0$ they
maintain the properties required in the corollary. 
\hfill$\Box$ 

\begin{cor} \label{kyft2}
There is a configuration of $k+3$ open (closed) sets 
covering $S^k$ none of which contains a pair of antipodal points, such
that no ${\mbf x}\in S^k$ is contained in more than $\lceil{{k+3}\over
2}\rceil$ of these sets and for every ${\mbf x}\in S^k$ the number of sets
that contain one of ${\mbf x}$ and $-{\mbf x}$ is at most $k+2$. 
\end{cor}

\proof For closed sets consider the sets $A_i$ in the proof of Corollary
\ref{kyftight} together with the closure of $D$. For open sets consider the
open $\varepsilon$-neighborhoods of these sets for suitably small
$\varepsilon>0$.  
\hfill$\Box$

Note that covering with $k+3$ sets is optimal in Corollary~\ref{kyft2}
if $k\ge 3$. By the Borsuk-Ulam Theorem (form (i)) fewer than $k+2$ open
(or closed) sets not containing antipodal pairs of points is not enough to
cover $S^k$. If we cover with $k+2$ sets (open or closed), then it gives rise
to a proper coloring of $B(k+1,\alpha)$ for large enough $\alpha$ in a natural
way. This coloring uses the optimal number $k+2$ of colors, therefore it has a
vertex with $k+1$ different colors in its neighborhood. A compactness argument 
establishes from this that there is a point in $S^k$ covered by $k+1$ sets.
A similar argument gives that $k+2$ in Corollary \ref{kyftight} is also
optimal if $k\ge 3$. 

\begin{cor} \label{qbound}
$${h\over 2}+1\leq Q(h)\leq {h\over 2}+2.$$
\end{cor}

\proof
The lower bound is implied by Ky Fan's theorem. The upper bound follows from  
Corollary \ref{kyft2}.
\hfill$\Box$

\medskip
Notice that for odd $h$ Corollary~\ref{qbound} gives the exact value
$Q(h)={h+3\over 2}$. For $h$ even we either have $Q(h)={h\over 2}+1$ or 
$Q(h)={h\over 2}+2$. It is trivial that $Q(0)=1$. In \cite{up} we will show
$Q(2)=3$. This was independently proved by Imre B\'ar\'any \cite b.
For $h>2$ even it remains open whether the lower or the upper bound of
Corollary~\ref{qbound} is exact.

\section{Circular colorings}\label{sec:circ}

In this section we show an application of the Zig-zag Theorem for the circular
chromatic number of graphs. This will result in the partial solution of a
conjecture by Johnson, Holroyd, and Stahl \cite{JHS} and in a partial answer
to a question of Hajiabolhassan and Zhu \cite{HZ} concerning the circular
chromatic number of Kneser graphs and Schrijver graphs, respectively. 
We also answer a question of Chang, Huang, and Zhu \cite{CHZ} concerning the
circular chromatic number of iterated Mycielskians of complete graphs. 

The circular chromatic number of a graph was introduced by Vince \cite{Vin}
under the name star chromatic number as follows.

\begin{defi} \label{defi:circ}
For positive integers $p$ and $q$ a coloring $c:V(G)\to [p]$ of a graph $G$ is
called a {\em $(p,q)$-coloring} if for all adjacent vertices $u$ and $v$ one
has $q\leq |c(u)-c(v)|\leq p-q$. The {\em circular chromatic number} of $G$ is
defined as
$$\chi_c(G)=\inf\left\{{p\over q}:\hbox{\rm\ there is a $(p,q)$-coloring of }
G\right\}.$$ 
\end{defi}

It is known that the above infimum is always attained for finite graphs. An
alternative
description of $\chi_c(G)$, explaining its name, is that it is the minimum
length of the perimeter of a circle on which we can represent the vertices of
$G$ by arcs of length $1$ in such a way that arcs belonging to adjacent
vertices do not overlap. For a proof of this equivalence and for an extensive
bibliography on the circular chromatic number we refer to Zhu's survey article
\cite{Zhu}.

It is known that for every graph $G$ one has 
$\chi(G)-1<\chi_c(G)\leq \chi(G)$. Thus
$\chi_c(G)$ determines the value of $\chi(G)$ while this is not true the other
way round. Therefore the circular chromatic number can be considered as a
refinement of the chromatic number. 
\medskip

Our main result on the circular chromatic number is
Theorem~\ref{thm:circ}. Here we restate the theorem with the explicit meaning
of being topologically $t$-chromatic.
\medskip

\noindent{\bf Theorem~\ref{thm:circ}} (restated) {\em For a finite graph $G$
we have $\chi_c(G)\ge\coind(B_0(G))+1$ if $\coind(B_0(G))$ is odd.
}
\medskip

\proof Let $t=\coind(B_0(G))+1$ be an even number and let $c$ be a
$(p,q)$-coloring of $G$. 
By the Zig-zag Theorem there is a $K_{{t\over 2},{t\over 2}}$ 
in $G$ which is completely multicolored by colors appearing in an alternating
manner in its two sides. Let these colors be
$c_1<c_2<\dots<c_t$. Since the vertex colored $c_i$ is adjacent to that
colored $c_{i+1}$, we have $c_{i+1}\ge c_i+q$ and
$c_t\ge c_1+(t-1)q$. Since $t$ is even, the vertices colored $c_1$ and $c_t$
are also adjacent, therefore we must have $c_t-c_1\leq p-q$. The last two
inequalities give $p/q\ge t$ as needed.
\hfill$\Box$

This result has been independently obtained by Meunier \cite{meunier} for
Schrijver graphs. 

\subsection{Circular chromatic number of even chromatic 
Kneser and Schrijver graphs} 

Johnson, Holroyd, and Stahl \cite{JHS} considered the circular chromatic
number of Kneser graphs and formulated the following conjecture. 
(See also as Conjecture 7.1 and Question 8.27 in \cite{Zhu}.)

\smallskip
\par
\noindent
{\bf Conjecture} (Johnson, Holroyd, Stahl \cite{JHS}):     
For any $n\ge 2k$ $$\chi_c(KG(n,k))=\chi(KG(n,k)).$$
\medskip
\par
\noindent
It is proven in \cite{JHS} that the above conjecture holds if $k=2$ or
$n=2k+1$ or $n=2k+2$. 

Lih and Liu \cite{LihLiu} investigated the circular chromatic number of
Schrijver graphs and proved that $\chi_c(SG(n,2))=n-2=\chi(SG(n,2))$ whenever
$n\neq 5$. (For $n=2k+1$ one always has $\chi_c(SG(2k+1,k))=2+{1\over k}$.) It
was conjectured in \cite{LihLiu} and proved in
\cite{HZ} that for every fixed $k$ there is a threshold $l(k)$ for which $n\ge
l(k)$ implies $\chi_c(SG(n,k))=\chi(SG(n,k))$. This clearly implies the
analogous statement for Kneser graphs, for which the explicit threshold 
$l(k)=2k^2(k-1)$ is given in \cite{HZ}. At the end of their paper \cite{HZ} 
Hajiabolhassan and Zhu ask what is the minimum $l(k)$ for which $n\ge l(k)$
implies $\chi_c(SG(n,k))=\chi(SG(n,k))$. We show that no such threshold
is needed if $n$ is even.

\begin{cor} \label{cor:JHS}
The Johnson-Holroyd-Stahl conjecture holds for every even $n$.
Moreover, if $n$ is even, then the stronger equality
$$\chi_c(SG(n,k))=\chi(SG(n,k))$$
also holds.
\end{cor}

\proof As $t$-chromatic Kneser graphs and Schrijver graphs are topologically
$t$-chromatic, Theorem~\ref{thm:circ} implies the statement of the corollary.
\hfill$\Box$

As mentioned above this result has been obtained independently by Meunier
\cite{meunier}. 

\smallskip

We show in Subsection \ref{oddsch} that for odd $n$ the situation is
different.

\subsection{Circular chromatic number of Mycielski graphs and Borsuk graphs}

The circular chromatic number of Mycielski graphs was also studied
extensively, cf.\ \cite{CHZ,Fan,HZMyc,Zhu}. Chang, Huang, and Zhu \cite {CHZ}
formulated the conjecture that $\chi_c(M^d(K_n))=\chi(M^d(K_n))=n+d$
whenever $n\ge d+2$. Here $M^d(G)$ denotes the $d$-fold iterated Mycielskian
of graph $G$, i.e., using the notation of Subsection~\ref{subsect:gmyc} we
have $M^d(G)=M_{\msbf r}^{(d)}(G)$ with ${\mbf r}=(2,\dots,2)$.    
The above conjecture was verified for the special cases $d=1,2$ in \cite{CHZ},
where 
it was also shown that $\chi_c(M^d(G))\leq \chi(M^d(G))-1/2$ if 
$\chi(G)=d+1$. 
A simpler proof for the above special cases of the conjecture 
was given (for $d=2$ with the extra condition $n\ge 5$) in \cite{Fan}. 
Recently Hajiabolhassan and Zhu \cite{HZMyc} proved that $n\ge 2^d+2$ implies
$\chi_c(M^d(K_n))=\chi(M^d(K_n))=n+d$. Our results show that
$\chi_c(M^d(K_n))=\chi(M^d(K_n))=n+d$ always  
holds if $n+d$ is even. This also answers the question of Chang, Huang, and
Zhu asking the value of $\chi_c(M^n(K_n))$ (Question 2 in
\cite{CHZ}). The stated equality is given by the following
immediate consequence of Theorem~\ref{thm:circ}. 

\begin{cor} \label{cor:Mycirc}
If $H(G)$ is homotopy equivalent to the sphere $S^h$, ${\mbf r}$ is a
vector of positive integers, and $h+d$ is even, then
$\chi_c(M_{\msbf r}^{(d)}(G))\ge d+h+2$.
\par
\noindent
In particular, $\chi_c(M_{\msbf r}^{(d)}(K_n))=n+d$ whenever $n+d$ is even.
\end{cor}
 
\proof
The condition on $G$ implies ${\rm coind}(H(M_{\msbf r}^{(d)}(G)))=h+d$ by 
Stiebitz's
result \cite{Stieb} (cf.\ the discussion and Proposition \ref{prop:Stieb}
in Subsection~\ref{subsect:gmyc}),
which further implies ${\rm coind}(B_0(M_{\msbf r}^{(d)}(G)))=h+d+1$. This gives
the conclusion by Theorem~\ref{thm:circ}. 

The second statement follows by the homotopy equivalence of $H(K_n)$ with
$S^{n-2}$ and the chromatic number of $M_{\msbf r}^{(d)}(K_n)$ being $n+d$.
\hfill$\Box$ 

\smallskip
\par
\noindent
The above mentioned conjecture of Chang, Huang, and Zhu for $n+d$ even is a
special case with ${\mbf r}=({2,2,\dots,2})$ and $n\ge d+2$. Since $n+n$ is
always even, the answer $\chi_c(M^n(K_n))=2n$ to their question also follows.  
\smallskip
\par
\noindent
Corollary~\ref{cor:Mycirc} also implies a recent result of Lam, Lin, Gu, and
Song \cite{LLGS} who proved that for the generalized Mycielskian of odd order
complete graphs $\chi_c(M_r(K_{2m-1}))=2m$. 
\medskip
\par
\noindent
Lam, Lin, Gu, and Song \cite{LLGS} also determined the circular chromatic
number of the generalized Mycielskian of even order complete graphs. They
proved $\chi_c(M_r(K_{2m}))=2m+1/(\lfloor(r-1)/m\rfloor+1)$.
This result can be used to bound the circular chromatic number of the
Borsuk graph $B(2s,\alpha)$ from above.

\begin{thm} \label{thm:Borcirc} For the Borsuk graph $B(n,\alpha)$ we have
\begin{description}
\item[(i)]
$\chi_c(B(n,\alpha))=n+1$
if $n$ is odd and $\alpha\ge\sqrt{2+2/n}$; 
\item[(ii)] $\chi_c(B(n,\alpha))\to n$ as $\alpha\to2$ if $n$ is even.
\end{description}
\end{thm}

\proof
The lower bound of part (i) immediately follows from Theorem \ref{thm:circ}
considering again the finite subgraph $G_P$ of $B(n,\alpha)$ defined in
\cite{LLgomb} 
and already mentioned in the proof of Lemma~\ref{gh}. 
The matching upper bound is provided by $\chi(B(n,\alpha))=n+1$ for the range
of $\alpha$ considered, see \cite{LLgomb}. 

For (ii) we have $\chi_c(B(n,\alpha))>\chi(B(n,\alpha))-1\ge n$. For an upper
bound we use that $\chi_c(M_r(K_n))\to n$ if $r$ goes to infinity by the
result of Lam, Lin, Gu, and Song \cite{LLGS} quoted above. By the result of
Stiebitz \cite{Stieb} and Lemma~\ref{gh} we have a graph homomorphism from
$B(n,\alpha)$ to $M_r(K_n)$ for any $r$ and large enough $\alpha$. As
$(p,q)$-colorings can be defined in terms of graph homomorphisms (see
\cite{BH}), we have $\chi_c(G)\le\chi_c(H)$ if there exists a graph 
homomorphism from $G$ to $H$. This finishes the proof of part (ii) of the
theorem. 
\hfill$\Box$

\medskip
\noindent
{\em Remark 7.} By Theorem~\ref{thm:Borcirc} (ii) we have a sequence of
$(p_i,q_i)$-colorings of the graphs $B(n,\alpha_i)$ where $n$ is even such
that $\alpha_i\to2$ and $p_i/q_i\to n$. By a direct construction we can show
that a single function $g:S^{n-1}\to C$ is enough. Here $C$ is a circle of
unit perimeter. We need
\begin{equation}\label{eqb}
\inf\{\hbox{\rm dist}_C(g({\mbf x}),g({\mbf
y})):\{{\mbf x},{\mbf y}\}\in E(B(n,\alpha))\}\to 1/n\hbox{ as 
$\alpha<2$ goes to }2.
\end{equation}
The distance ${\rm dist}_C(\cdot,\cdot)$ is measured
along the  cycle $C$. Clearly, if $p/q>n$ and we split $C$ into $p$ arcs
$a_1,\ldots,a_p$ of equal length and color the point ${\mbf x}$ with $i$ if
$g({\mbf x})\in a_i$, then this is a $(p,q)$-coloring of $B(n,\alpha)$ for
$\alpha$ close enough to $2$.

For $n=2$ any $\2$-map $g:S^1\to C$ satisfies Equation~(\ref{eqb}). Let
$n>2$. The map $g$ to be constructed must not be continuous by the Borsuk-Ulam
theorem. Let us choose a set $H$ of $n-1$ equidistant points in $C$ and for
$b\in C$ let $T(b)$ denote the unique set of $n/2$ equidistant points in $C$
containing $b$.

We consider $S^{n-1}$ as the {\em join} of the sphere $S^{n-3}$ and the
circle $S^1$. All points in $S^{n-1}$ are now either in $S^{n-3}$, or in
$S^1$, or in the interval connecting a point in $S^{n-3}$ to a point in $S^1$.
We define $g$ on $S^{n-3}$ such that it takes values only from $H$ and it is a
proper coloring of $B(n-2,1.9)$. We define $g$ on $S^1$ such that if $\mbf y$
goes a full circle around $S^1$ with uniform velocity, then its image $g({\mbf
y})$ covers an arc of length $2/n$ of $C$
and it also moves with uniform velocity. Notice
that although $g$ is not continuous on $S^1$, the set $T(g({\mbf y}))$
depends on ${\mbf y}\in S^1$ in a continuous manner. Also note that for a point
${\mbf x}\in S^1$ the images $g({\mbf x})$ and $g(-{\mbf x})$ are $1/n$ apart on
$C$ and $T(g({\mbf x}))\cup T(g(-{\mbf x}))$ is a set of $n$ equidistant
points.

Let ${\mbf x}\in S^{n-3}$ and ${\mbf y}\in S^1$. Assume that a point $\mbf z$
moves with uniform velocity from $\mbf x$ to $\mbf y$ along the interval
connecting them. We define $g$ on this interval such that $g({\mbf z})$ moves
with uniform velocity along $C$ covering an arc of length at most $1/n$ from
$g({\mbf x})$ to a point in $T(g({\mbf y}))$. The choice of the point in
$T(g({\mbf y}))$ is uniquely determined unless $g({\mbf x})\in T(g(-{\mbf
x}))$. In the latter case we make an arbitrary choice of the two possible
points for the destination of the image $g({\mbf z})$.

It is not hard to prove that the function $g$ defined above satisfies
Equation~(\ref{eqb}).
\hfill$\Diamond$

\subsection{Circular chromatic number of odd chromatic Schrijver graphs}
\label{oddsch}

In this subsection we show that the parity condition on $\chi(SG(n,k))$ in 
Corollary \ref{cor:JHS} is relevant, for odd chromatic Schrijver graphs the
circular chromatic number can be arbitrarily close to its lower bound. 

\begin{thm} \label{thm:oddsch}
For every $\varepsilon>0$ and every odd $t\ge3$ if $n\ge t^3/\varepsilon$ and
$t=n-2k+2$, then 
$$1-\varepsilon<\chi(SG(n,k))-\chi_c(SG(n,k))<1.$$
\end{thm}

The second inequality is well-known and holds for any graph. We included
it only for completeness. To prove the first inequality we need some
preparation. We remark that the bound on $n$ in the theorem is not best
possible. Our method proves $\chi(SG(n,k))-\chi_c(SG(n,k))\ge1-1/i$ if $i$ is
a positive integer and $n\ge6(i-1){t\choose3}+t$.

First we extend our notion of wide coloring. 

\begin{defi}
For a positive integer $s$ we call a vertex coloring of a graph $s$-wide if
the two end vertices of any walk of length $2s-1$ receive different colors. 
\end{defi}

Our original wide colorings are $3$-wide, while $1$-wide simply means
proper. Gy\'arf\'as, Jensen, and Stiebitz \cite{GyJS} investigated $s$-wide
colorings (in different terms) and mention (referring to a referee in the
$s>2$ case) the existence of homomorphism universal graphs for $s$-wide
colorability with $t$ colors. We give a somewhat different family of such
universal graphs. In the $s=2$ case the color-criticality of
the given universal graph is proven in \cite{GyJS} implying its minimality
among graphs admitting $2$-wide $t$-colorings. 
Later in Subsection~\ref{critical} we generalize this result 
showing that the members of our family are color-critical for every $s$. 
Thus they must be minimal and therefore isomorphic to a retract of the
corresponding graphs given in \cite{GyJS}. 

\begin{defi} \label{Wst}
Let $H_s$ be the path on the vertices $0,1,2,\dots,s$ ($i$ and $i-1$ connected
for $1\le i\le s$) with a loop at $s$.
We define $W(s,t)$ to be the graph with 
$$V(W(s,t))=\{(x_1\dots x_t): \forall i\ x_i\in\{0,1,\dots,s\},\exists!i\
x_i=0,\ \exists j\ x_j=1\},$$
$$E(W(s,t))=\{\{x_1\dots x_t,y_1\dots y_t\}: \forall i\ \{x_i,y_i\}\in
E(H_s)\}.$$ 
\end{defi}

Note that $W(s,t)$ is an induced subgraph of the direct power $H_s^t$.

\begin{prop} \label{Wuni}
A graph $G$ admits an $s$-wide coloring with $t$ colors if and only if there
is a homomorphism from $G$ to $W(s,t)$.
\end{prop}

\proof
For the if part color vertex ${\mbf x}=x_1\dots x_t$ of $W(s,t)$ with $c({\mbf
x})=i$ if 
$x_i=0$. Any walk between two vertices colored $i$ either has even length or
contains two vertices ${\mbf y}$ and ${\mbf z}$ with $y_i=z_i=s$. These ${\mbf y}$
and ${\mbf z}$ are both at least at distance $s$ apart from both ends of the
walk, thus our 
coloring of $W(s,t)$ with $t$ colors is $s$-wide. Any graph admitting a
homomorphism $\varphi$ to $W(s,t)$ is $s$-widely colored with $t$ colors by
$c_G(v):=c(\varphi(v))$.  

For the only if part assume $c$ is an $s$-wide $t$-coloring of $G$ with
colors $1,\dots,t$. Let $\varphi(v)$ be an arbitrary vertex of $W(s,t)$ if $v$
is an isolated vertex of $G$. For a non-isolated vertex $v$ of $G$ let
$\varphi(v)={\mbf x}=x_1\dots x_t$ with $x_i=\min(s,d_i(v))$, where $d_i(v)$ is
the distance of color class $i$ from $v$. It is clear that $x_i=0$ for
$i=c(v)$ and for no other $i$, while $x_i=1$ for the colors of the neighbors
of $v$ in $G$. Thus the
image of $\varphi$ is indeed in $V(W(s,t))$. It takes an easy checking 
that $\varphi$ is a homomorphism. 
\hfill$\Box$

\medskip

The following lemma is a straightforward extension of the argument given in
the proof of Theorem \ref{thm:upb}. 

\begin{lem} \label{SGswide}
If $n\ge (2s-2)t^2-(4s-5)t$ then $SG(n,k)$ admits an $s$-wide $t$-coloring.  
\end{lem}

\proof 
We use the notation introduced in the proof of Theorem \ref{thm:upb}. 

Let $n\ge t(2(s-1)(t-2)+1)$ as in the statement and let $c_0$ be the coloring
defined in the mentioned proof. The lower bound on $n$ now allows to assume
that $|C_i|\ge(s-1)(t-2)+1$. We show that $c_0$ is $s$-wide.

Consider a walk $x_0x_1\dots x_{2s-1}$ of length $(2s-1)$
in $SG(n,k)$ and let $i=c_0(x_0)$. Then $C_i\subseteq x_0$. By Lemma
\ref{partav} $|x_0\setminus x_{2s-2}|\leq(s-1)(t-2)<|C_i|$. Thus $x_{2s-2}$ is
not disjoint from $C_i$. As $x_{2s-1}$ is disjoint from $x_{2s-2}$, it does
not contain $C_i$ and thus its color is not $i$.
\hfill$\Box$

\begin{lem} \label{WM}
$W(s,t)$ admits a homomorphism to $M_s(K_{t-1})$. 
\end{lem}
 
\proof
Recall our notation for the (iterated) generalized Mycielskians from Subsection
\ref{subsect:gmyc}.  

We define the following mapping from $V(W(s,t))$ to $V(M_s(K_{t-1}))$. 

$$\varphi(x_1\dots x_t):=\left\{\begin{array}{lll}(s-x_t,i)&&\hbox{if }
x_t\neq x_i=0\\(s,*)&&\hbox{if }x_t=0.\end{array}\right.$$

One can easily check that $\varphi$ is indeed a homomorphism.
\hfill$\Box$

\medskip
\par
\noindent
{\bf Proof of Theorem \ref{thm:oddsch}.}
By Lemma \ref{SGswide}, if $n\ge(2s-2)t^2-(4s-5)t$, then $SG(n,k)$ has an
$s$-wide $t$-coloring, thus by
Proposition \ref{Wuni} it admits a homomorphism to $W(s,t)$. Composing this
with the homomorphism given by Lemma \ref{WM} we conclude that $SG(n,k)$ admits
a homomorphism to $M_s(K_{t-1})$, implying
$\chi_c(SG(n,k))\leq\chi_c(M_s(K_{t-1}))$. 

We continue by using Lam, Lin, Gu, and Song's result \cite{LLGS}, who
proved, as already quoted in the previous subsection, that  
$\chi_c(M_s(K_{t-1}))=t-1+{1\over\left\lfloor2s-2\over t-1\right\rfloor+1}$ if
$t$ is odd. Thus, for odd $t$ and $i>0$ integer we choose $s=(t-1)(i-1)/2+1$
and $\chi(SG(n,k))-\chi_c(SG(n,k))=t-\chi_c(SG(n,k))\ge1-1/i$ follows from the
$n\ge6(i-1){t\choose3}+t$ bound.

To get the form of the statement claimed in the theorem we choose
$i=\lfloor1/\varepsilon\rfloor+1$.
\hfill$\Box$

\section{Further remarks}

\subsection{Color-criticality of $W(s,t)$}\label{critical}

In this subsection we prove the edge color-criticality of the graphs $W(s,t)$
introduced in the previous section. This generalizes Theorem 2.3 in
\cite{GyJS}, see Remark 8 after the proof. 

\begin{thm}
For every integer $s\ge 1$ and $t\ge 2$ the graph $W(s,t)$ has chromatic
number $t$, but deleting any of its edges the resulting graph is
$(t-1)$-chromatic.  
\end{thm}

\proof
$\chi(W(s,t))\ge t$ follows from the fact that some $t$-chromatic Schrijver 
graphs admit
a homomorphism to $W(s,t)$ which is implied by Lemma \ref{SGswide} and
Proposition \ref{Wuni}. The coloring giving vertex ${\mbf x}=x_1\dots x_t$
of $W(s,t)$ color $i$ iff $x_i=0$ is proper proving $\chi(W(s,t))\leq t$.

We prove edge-criticality by induction on $t$. 
For $t=2$ the statement is trivial as $W(s,t)$ is isomorphic to $K_2$. 
Assume that $t\ge 3$ and edge-criticality holds for $t-1$. Let $\{x_1\dots
x_t,y_1\dots y_t\}$ be  
an edge of $W(s,t)$ and $W'$ be the graph remaining after removal of this
edge. We need to give a proper $(t-1)$-coloring $c$ of $W'$.

Let $i$ and $j$ be the coordinates for which $x_i=y_j=0$. We have $x_j=y_i=1$,
in particular, $i\neq j$.  
Let $r$ be a coordinate different from both $i$ and $j$. We may assume without
loss of generality that $r=1$, and also that $y_1\ge x_1$. Coordinates $i$ and
$j$ make sure that $x_2x_3\dots x_t$ and $y_2y_3\dots y_t$ are vertices of
$W(s,t-1)$, and in fact, they are connected by an edge $e$.

A proper $(t-2)$-coloring of the graph $W(s,t-1)\setminus e$ exists by the
induction hypothesis. Let $c_0$ be such a coloring. Let $\alpha$ be a color of
$c_0$ and $\beta$ a color that does not appear in $c_0$. We define the
coloring $c$ of $W'$ as follows:
$$c(z_1z_2\dots z_t)=\left\{\begin{array}{lll}
\alpha&&\hbox{if }z_1<x_1,\ x_1-z_1\hbox{ is even}\\
\beta&&\hbox{if }z_1<x_1,\ x_1-z_1\hbox{ is odd}\\
\alpha&&\hbox{if }z_1=x_1=1,\ z_i\ne1\hbox{ for }i>1\\
\beta&&\hbox{if }z_1>x_1,\ z_i=x_i\hbox{ for }i>1\\
c_0(z_2z_3\dots z_t)&&\hbox{otherwise.}
\end{array}\right.$$

It takes a straightforward case analysis to check that $c$ is a proper
$(t-1)$-coloring of $W'$.   
\hfill$\Box$

\medskip
\par
\noindent
{\em Remark 8.} Gy\'arf\'as, Jensen, and Stiebitz \cite{GyJS} 
proved the $s=2$ version of
the previous theorem using a homomorphism from their universal graph with
parameter $t$ to a generalized Mycielskian of the same type of graph with
parameter $t-1$. In fact, our proof is a direct generalization of theirs using
very similar ideas. Behind the coloring we gave is the recognition of a
homomorphism from $W(s,t)$ to $M_{3s-2}(W(s,t-1))$. 
\hfill$\Diamond$

\subsection{Hadwiger's conjecture and the Zig-zag theorem}

Hadwiger's conjecture, one of the most famous open problems in graph theory,
states that if a graph $G$ contains no $K_{r+1}$ minor, then $\chi(G)\leq r$. 
For detailed information on the history and status of this conjecture we refer
to Toft's survey \cite{Toft}. We only mention that even $\chi(G)=O(r)$ is not
known to be implied by the hypothesis for general $r$. 

As a fractional and linear approximation version, Reed and Seymour \cite{RS}
proved that if $G$ has no $K_{r+1}$ minor then $\chi_f(G)\leq 2r$. This means
that graphs with $\chi_f(G)$ and $\chi(G)$ appropriately close and not
containing a $K_{r+1}$ minor satisfy $\chi(G)=O(r)$. 

We know that the main examples of graphs in \cite{SchU} for
$\chi_f(G)<<\chi(G)$ (Kneser graphs, Mycielski graphs) satisfy the hypothesis
of the Zig-zag theorem, therefore their $t$-chromatic versions must contain 
$K_{\lceil{t\over 2}\rceil,\lfloor{t\over 2}\rfloor}$ subgraphs. 
(We mention that for strongly topologically $t$-chromatic graphs this
consequence, in fact, the containment of $K_{a,b}$ for every $a,b$ satisfying
$a+b=t$, was proven by Csorba, Lange, Schurr, and Wassmer \cite{CsLSW}.)  
However, a $K_{\lceil{t\over 2}\rceil,\lfloor{t\over 2}\rfloor}$
subgraph contains a $K_{\lfloor{t\over 2}\rfloor +1}$ minor (just take a
matching of size $\lfloor{{t-2}\over 2}\rfloor$ plus one
point from each side of the bipartite graph) proving the
following statement which shows that the same kind of approximation is valid
for these graphs, too. 

\begin{cor} \label{Hadw}
If a topologically $t$-chromatic graph contains no $K_{r+1}$ minor, then
$t<2r.$
\end{cor} 
\hfill$\Box$

\bigskip
\par
\noindent
{\bf Acknowledgments:}
We thank Imre B\'ar\'any, P\'eter Csorba, G\'abor Elek, L\'aszl\'o Feh\'er,
L\'aszl\'o Lov\'asz, Ji\v{r}\'{\i} Matou\v{s}ek, and G\'abor Moussong for many
fruitful conversations that helped us to better understand the topological
concepts used in this paper.

\end{document}